\documentclass[preprint,12pt,authoryear]{elsarticle}

\usepackage{amssymb}
\usepackage{amsmath}
\usepackage[english]{babel}
\usepackage{amsthm}   
\usepackage{mathabx}  
\usepackage{graphicx} 
\usepackage{fancyhdr} 
\usepackage{setspace} 
\usepackage{lineno}	  
\usepackage{lastpage} 
\usepackage{booktabs} 
\usepackage{rotating} 
\usepackage{multirow}	
\usepackage{caption,subcaption}  
\usepackage[round]{natbib}   

\usepackage{comment}
\usepackage{siunitx}
    \DeclareSIUnit{\gramdryweight}{\ensuremath{g_{CDW}}}
    \DeclareSIUnit{\year}{\ensuremath{a}}
    \DeclareSIUnit{\bar}{\ensuremath{bar}}
    \DeclareSIUnit{\USD}{\ensuremath{USD}}
\usepackage{longtable}
\usepackage{tabularx}
\usepackage[version=4]{mhchem}
\usepackage{diagbox}
\usepackage{url}

\begin{document}

\begin{frontmatter}

\title{Simultaneous Design of Microbe and Bioreactor}
\author[SVT]{Anita L. Ziegler}
\author[SVT]{Marc-Daniel Stumm}
\author[SVT]{Tim Prömper}
\author[BioVT]{Thomas Steimann}
\author[BioVT]{J{\o}rgen Magnus}
\author[JARA,SVT,ICE-1]{Alexander Mitsos \corref{cor1}} 

\ead{amitsos@alum.mit.edu}
\cortext[cor1]{Corresponding author}

\affiliation[SVT]{organization={Process Systems Engineering (AVT.SVT), RWTH Aachen University},
	city={Aachen},
	postcode={52074}, 
	country={Germany}}
\affiliation[BioVT]{organization={Biochemical Engineering (AVT.BioVT), RWTH Aachen University},
	city={Aachen},
	postcode={52074}, 
	country={Germany}}
\affiliation[JARA]{organization={JARA-ENGERGY},
	city={Aachen},
	postcode={52056}, 
	country={Germany}}
\affiliation[ICE-1]{organization={Institute of Climate and Energy Systems, Energy Systems Engineering (ICE-1), Forschungszentrum Jülich GmbH},
	city={Jülich},
	postcode={52425}, 
	country={Germany}}

\begin{abstract} 
When developing a biotechnological process, the microorganism is first designed, e.g., using metabolic engineering.
Then, the optimum fermentation parameters are determined on a laboratory scale, and lastly, they are transferred to the bioreactor scale. 
However, this step-by-step approach is costly and time-consuming and may result in suboptimal configurations. 
Herein, we present the bilevel optimization formulation \textit{SimulKnockReactor}, which connects bioreactor design with microbial strain design, an extension of our previous formulation, SimulKnock (Ziegler et al., 2024, AIChE J.). 
At the upper (bioreactor) level, we minimize investment and operation costs for agitation, aeration, and pH control by determining the size and operating conditions of a continuous stirred-tank reactor---without selecting specific devices like the stirrer type. 
The lower (cellular) level is based on flux balance analysis and implements optimal reaction knockouts predicted by the upper level.  
Our results with a core and a genome-scale metabolic model of \textit{Escherichia coli} show that the substrate is the largest cost factor. 
Our simultaneous approach outperforms a sequential approach using OptKnock. 
Namely, the knockouts proposed by OptKnock cannot guarantee the required production capacity in all cases considered. 
In the case that both approaches deliver feasible results, the total annual costs are the same or lower with SimulKnockReactor, highlighting the advantage of combining cellular and bioreactor levels.
This work is a further step towards a fully integrated bioprocess design. 

\end{abstract}
\begin{highlights}
    \item SimulKnockReactor connects the micro-level with the macro-level of a bioprocess.
    \item SimulKnockReactor proposes optimal gene deletions and minimized bioreactor costs. 
\item CAPEX and OPEX are considered for reactor, oxygen supply, pH control, and cooling.  
\item SimulKnockReactor achieves lower bioreactor costs compared to sequential design.
\item The implementation of SimulKnockReactor is openly available.
\end{highlights}

\begin{graphicalabstract}
    \includegraphics[width=\linewidth]{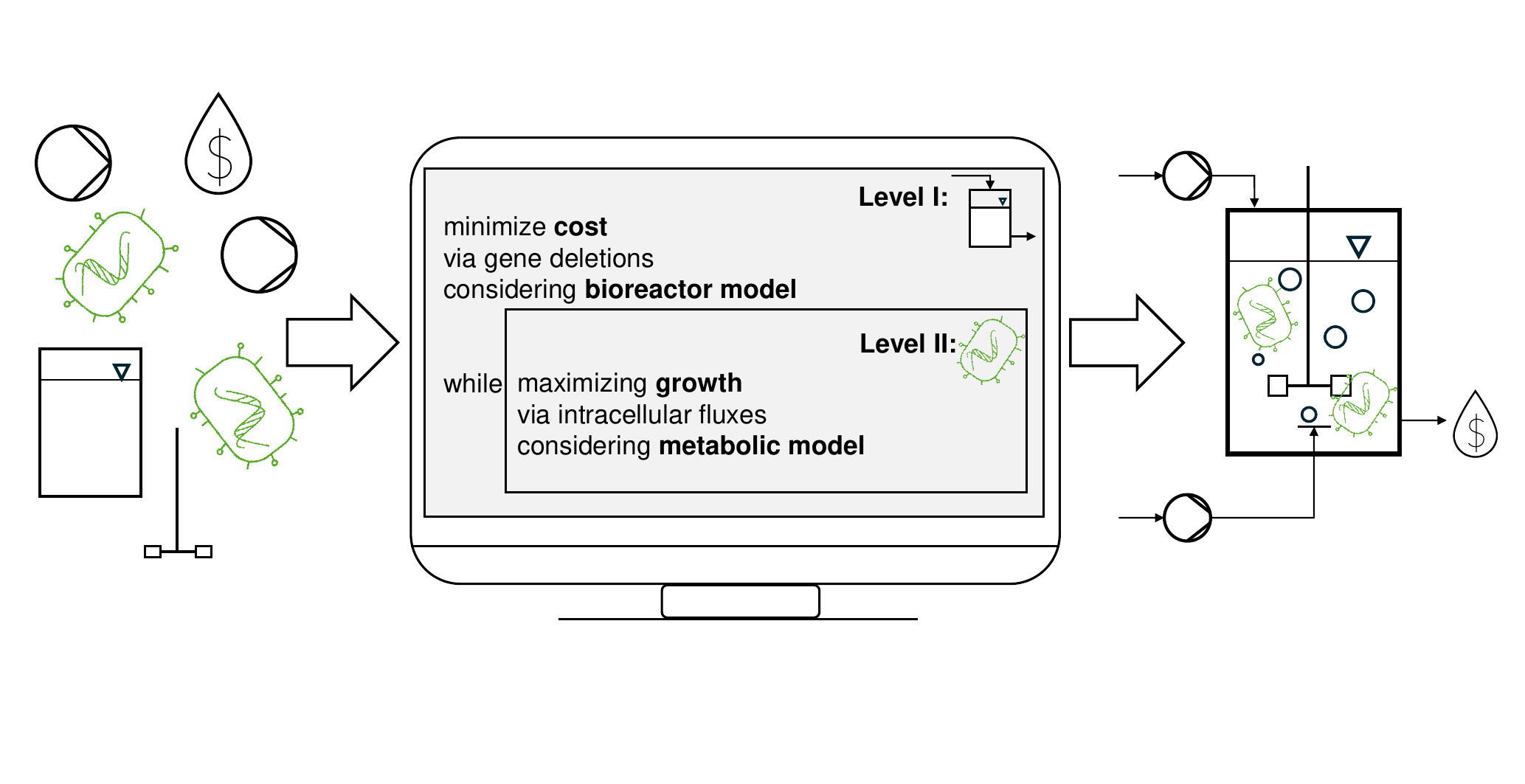}
\end{graphicalabstract}
\end{frontmatter}
\clearpage

\section{Introduction}
Model-based simulation and optimization are established methods for the design and scale-up of fermentation processes.
When modeling a fermentation process, different scales and their respective variables and parameters come into play, e.g., the cellular level, the fermentation reaction level, and the bioreactor level \citep{Pinto.2016, Villadsen.2011}. 

At the cellular level, the intracellular fluxes can be determined via metabolic modeling, typically via flux balance analysis (FBA) \citep{Orth.2010, Watson.1986, Savinell.1992}, which is based on steady-state mass balances of the intracellular metabolites. 
The metabolites and their corresponding reaction stoichiometry are contained in the genome-scale metabolic model (GEM) of the microorganism. 
With evolving available information in the GEM, different aspects have been included in metabolic modeling formulations, e.g., the intracellular pH \citep{Du.2019}, thermodynamics \citep{Niebel.2019}, and regulatory aspects, reviewed by \citet{Gombert.2000}.
First introduced by \citet{Burgard.2003}, numerous strain optimization formulations have been developed, that allow for predicting optimal genetic modifications for improved product yield.

At the fermentation reaction level, the concentrations of the representative extracellular compounds are described by mass balances together with whole-cell growth or enzyme kinetics. Often energy balances are also required.
Momentum balances can be set up \citep{Liu.2020} but typically lumped models are used.
In the most basic case, the representative compounds are biomass, substrate, and product. 
Depending on the fermentation mode, e.g., batch, fed-batch, or continuous, the mass balances are differential, algebraic, or differential-algebraic equations \citep{Chmiel.2018, Bailey.1986}. 
Exemplary fermentation parameters are temperature, pH, and nutrient composition \citep{Rao.2006, PortodeSouzaVandenberghe.2022}. 
\citet{Gordeeva.2013} presented one of the few examples of optimal fermentation design, namely for continuous fermentation including a recycling stream.
Another example was presented by \citet{Sinner.2019}, who proposed a bi-objective optimization formulation for maximal space-time yield and minimal sugar concentration in the effluent of continuous fermentation. 
To ensure and boost the performance of the fermentation, however, it is fruitful to connect the cellular level with the process levels \citep{Richelle.2020, Olsson.2022}.
\citet{Mahadevan.2002} were the first to connect the fermentation reaction level with the cellular level in their dynamic flux balance analysis. 
\citet{Chang.2016}, \citet{Jabarivelisdeh.2020}, and \citet{Bhonsale.2022} extended this connection to model-predictive control of the fermentation, similar to \citet{Oliveira.2021} and \citet{EspinelRios.2024}, who used machine learning techniques for substituting the cellular level. 
\citet{Zhuang.2013} considered yield, titer, and space-time yield for batch and fed-batch fermentation with product inhibition in their iterative strain design approach. 
\citet{Jabarivelisdeh.2018} also maximized the space-time yield, namely in a dynamic bilevel optimization program, switching on genetic modifications that were found in the preceding step.
In our previous work \citep{Ziegler.2024}, we combined strain optimization with optimal fermentation design.
In our bilevel optimization formulation \textit{SimulKnock}, the upper (fermentation) level includes the mass balances of continuous fermentation, while the lower (cellular) level is derived from FBA. 
SimulKnock was solved for two case studies with the GEM \textit{i}ML1515 \citep{Monk.2017} of \textit{E. coli} embedded. 
The results indicated that different knockouts were found and that the space-time yield could be significantly elevated compared to sequential optimization. 
This sequential optimization represents the common practice of first optimizing the microbial strain and then, subsequently, optimizing the fermentation. 
The optimization formulation SimulKnock, however, focuses on the cellular and the fermentation reaction level and does not give any information about the bioreactor level, including sizing of the bioreactor and cost. 

At the bioreactor level, in principle, the selection of reactor types, economic evaluation, and scale-up aspects become important \citep{Posten.2018, Clarke.2013}, to give a comprehensive list.
The selection of reactor type includes stirred tank reactors, pneumatically driven reactors, and other types \citep{Posten.2018}. 
Economic evaluation includes operating costs, material costs, separation costs, and product value, among others \citep{Liu.2020, Doran.1995}.
So far, the idea behind scale-up is to find suitable fermentation conditions at the laboratory scale and transfer these to the industrial bioreactor scale \citep{PortodeSouzaVandenberghe.2022,Clarke.2013}. 
Based on the parameters of the fermentation reaction, e.g., temperature and pH, the scale-up aspects are mechanical agitation, aeration, heat transfer, power consumption, rheology, and biotransformation \citep{Rao.2006, Villadsen.2011, Mandenius.2016}. 
Computational fluid dynamics have been used to model the spatial distribution of compounds in the bioreactor \citep{Kelly.2008, Rathore.2016}.
\citet{Cardoso.2020} and \citet{Knoll.2005} estimated bioreactor costs for different cultivation conditions. 
In optimal control, \citet{Hebing.2020} considered the oxygen transfer rate, and \citet{Stosch.2016} also included data on temperature, pH, and agitation rate measurement in their hybrid modeling approach. 
With fixed reactor volume, \citet{Powell.2008} chose a bioreactor setup and reduced cost in an iterative process.
\citet{Nath.2011} reviewed factorial design methods for improved hydrogen production. 
Interestingly, we did not find studies that used numerical optimization for optimal bioreactor design. 

In summary, the field of model-based optimal bioreactor design and cost estimation of bioreactors is underrepresented in biochemical engineering and, more specifically, the bioreactor level and the cellular level are yet to be connected. 
This connection is crucial in our eyes in order to design economically viable processes. 
An organism that reaches a high product yield and even a high space-time yield might not prove economically reasonable in an industrial bioreactor because other aspects like substrate cost and operation cost become decisive at the bioreactor level. 
In another case, the optimal fermentation conditions of the laboratory scale prove to not be realizable in the bioreactor. 
With the usual approach of optimizing the microbe, the fermentation conditions, and the bioreactor one after the other, it is necessary in these cases to go back and redevelop the microbial strain and the fermentation conditions, which is costly and time-consuming. 
By connecting strain design with fermentation and bioreactor design, we can circumvent this loop and simultaneously optimize the strain for the reactor and vice versa. 

Herein, we present \textit{SimulKnockReactor}, a methodology for bioreactor cost optimization considering gene deletion strategies. 
It combines bioreactor design with strain design in a bilevel optimization formulation and is an extension of our recent contribution SimulKnock \citep{Ziegler.2024}.
Within SimulKnockReactor, the lower level calculates the intracellular fluxes based on FBA, analogous to SimulKnock. 
The genetic modification, which we focus on, is gene deletion.
The upper level of SimulKnockReactor minimizes the cost of a bioreactor based on capital expenditures (CAPEX) and operational expenses (OPEX)---in contrast to SimulKnock, where the space-time yield was maximized.
The bioreactor is modeled as a set of mass and energy balances, combined with cost correlations. 
We estimate the cost for a stirred-tank reactor, for substrates, oxygen supply, heat transfer, and pH regulation.
We do not select specific devices of the bioreactor and the surrounding equipment; for example, we do not distinguish between a cooling jacket and a cooling coil or different stirrer types. 
Despite this fact, we still use the term bioreactor design throughout this work, since we do calculate the mass and energy flow rates of the fermentation, and the required mass flow rates of pH control agents, air, and cooling agent. 
Moreover, based on these flows, the required power for pumping, stirring, and compression, as well as the sizing of the reactor are determined. 
Information about oxygen and substrate consumption, biomass and product generation are retrieved from the lower (cellular) level.
This approach is expected to improve the accuracy compared to common correlations and allows to examine a specific organism without experimental setup. 
As exemplary case studies, we apply SimulKnockReactor to a core network and a GEM of \textit{E. coli} for aerobic production of formate, acetate, and succinate and compare the reaction deletions with SimulKnock results.
As discussed in the outlook, the reactor design could be extended to account for other aspects.

The contribution of this work is to provide a methodology for the initial assessment of bioreactor costs, reactor sizing, and operation, and simultaneously, the prediction of optimal gene deletions. 
The outstanding feature of the methodology is to connect the micro- and the macro-level of a bioprocess, or, in other words, the field of microbial strain design with bioreactor design. 
Our methodology is applicable to diverse microorganisms and target chemicals. 
Thus, this work is the next step towards a fully integrated bioprocess design, considering all levels of bioprocess engineering.

\section[Method]{Method: Simultaneous bioreactor design and microbial strain design}
\begin{figure}[htb]
    \centering
    \includegraphics[trim={0 4cm 0 3cm}, clip, width=\textwidth]{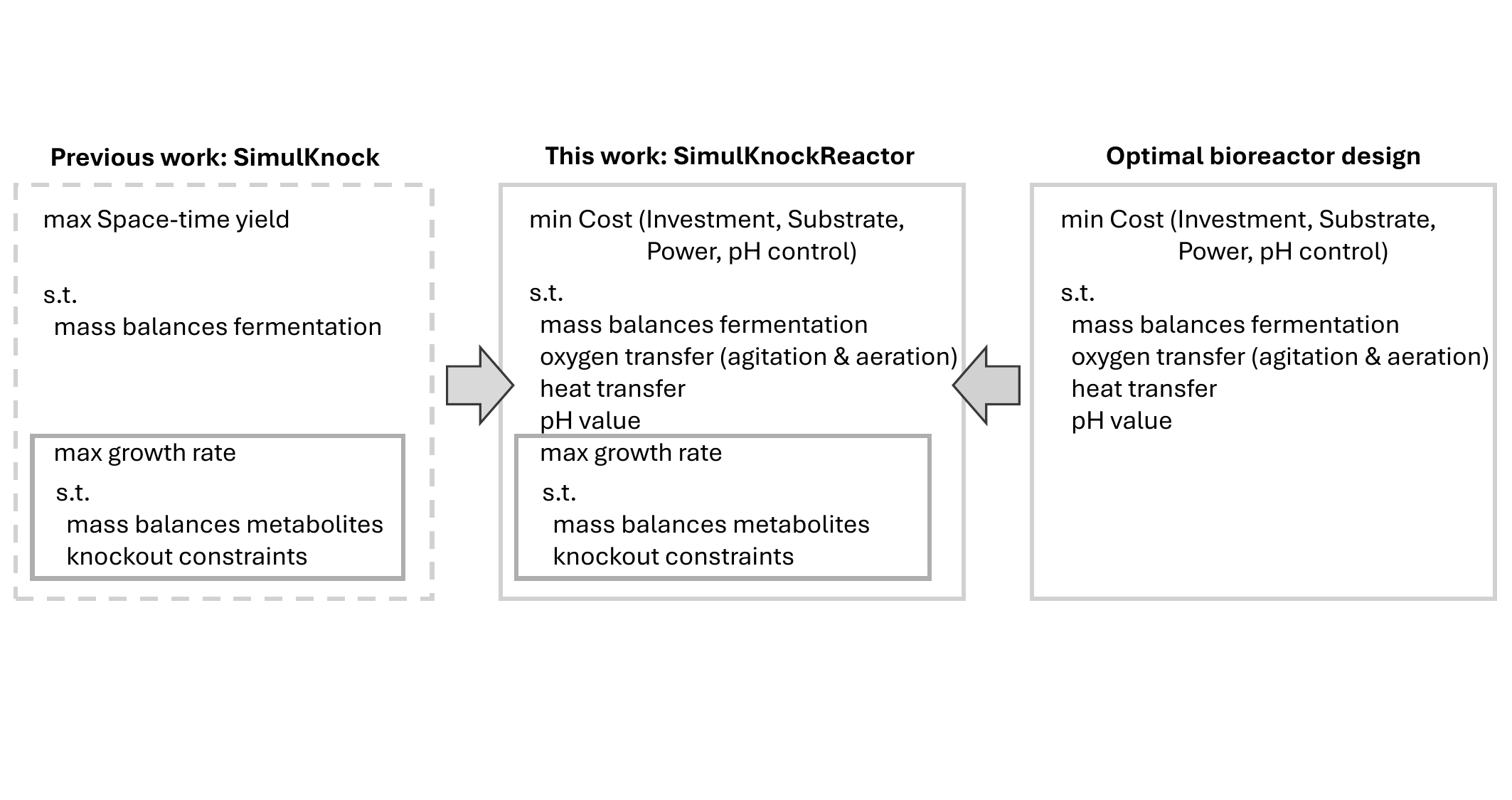}
    \caption{Overview of the elements of SimulKnock \citep{Ziegler.2024} and optimal bioreactor design that were combined to create the proposed formulation SimulKnockReactor in this work. The bilevel structure is adopted from SimulKnock:  the lower level  remains unchanged; the upper level is extended to include optimal bioreactor design. Note that optimal bioreactor design is not standard in the literature, however, the incorporated equations can largely be found in standard bioreactor and bioreaction literature \citep{Doran.1995, Villadsen.2011, Chmiel.2018, Liu.2020}}
\label{fig:overview_SimulKnock_Ext}
\end{figure}

Our aim was to minimize the cost of a bioreactor based on process variables and reaction knockouts concerning the microbial fluxes. 
The cost of the bioreactor consists of CAPEX and OPEX, of which we considered the aspects of aeration, agitation, pH control, and cooling. 
To calculate these costs, a model of the bioreactor and the surrounding equipment is needed which includes the compounds of the fermentation.
The fermentation, in turn, is coupled with the microbial metabolism by the uptake and secretion reactions of the organism. 
Our previous formulation SimulKnock \citep{Ziegler.2024} already represents a framework where the fermentation equations are coupled with the metabolic fluxes.  
In bilevel optimization, one optimization program which is called a lower-level program is embedded in another optimization program, which is called an upper-level program \citep{Dempe.2002}. 
In this work, we build upon SimulKnock by extending the model to allow for optimal bioreactor design.
Figure~\ref{fig:overview_SimulKnock_Ext} depicts the elements of SimulKnock and optimal bioreactor design and how we combined them to form SimulKnockReactor.  

We retained the bilevel structure and the lower level of SimulKnock (see Figure~\ref{fig:overview_SimulKnock_Ext}, left side) and extended the upper level of SimulKnock with a formulation of optimal bioreactor design. 
As a result, the objective of optimization changes from maximizing the space-time yield to minimizing costs.
Moreover, we included more fermentation compounds, namely oxygen, hydrogen, and strong ions. 

In contrast to a formulation of optimal bioreactor design (see Figure~\ref{fig:overview_SimulKnock_Ext}, right side), which is single-level, we used the lower-level variables instead of applying commonly used correlations.
This replacement is possible because the metabolic network includes exchange reactions that connect the metabolism to the environment, i.e., the fermentation culture.
These uptake and excretion reactions are also part of the bioreactor equations. 
By using lower-level variables in the upper level, we connected the bioreactor level with the cellular level.

In the following, the previous formulation of SimulKnock is described. 
Then, the objective and constraints for bioreactor design are introduced and merged with SimulKnock in the next step.
The implementation and solution of the resulting bilevel formulation with global and local optimization techniques are shown in the last section of this chapter. 

\subsection[SimulKnock]{Previous work: Simultaneous design of fermentation and microbe}

The bilevel optimization formulation SimulKnock \citep{Ziegler.2024} simultaneously optimizes a continuous fermentation process and a microbial strain for growth-coupled production.  
The upper level of SimulKnock represents the fermentation level, while the lower level represents the cellular level.
The connection between the upper and the lower levels is achieved via the implementation of kinetics and via the replacement of common correlations with cellular variables. 
We provide an overview of the SimulKnock formulation in the following and refer the interested reader to our previous publication \citep{Ziegler.2024} for details and assumptions. 
The mathematical formulation of SimulKnock with the embedded Monod kinetics reads

\begin{align}
    \underset{\substack{\boldsymbol{y},c_{S,Feed}, c_P, \\ c_{bio},  v_S', v_{bio}', v_P'}}\max \quad &c_{P}\cdot v_{bio}'  \quad \rightarrow\text{maximize space-time yield} \label{eq:max_STY}\\
    \text{s.t. }                                                                                   \kappa \geq& \sum_{i=1}^{r} (1-y_i)   \rightarrow\text{maximum number of knockouts}\label{eq:number_KO}\\
                                            0=&-v_{S}'\cdot M_S \cdot c_{bio} +v_{bio}'\cdot(c_{S,Feed} \nonumber\\
                                             &-K_{S} \frac{v_{bio}'}{v_{bio}^{max}-v_{bio}'}) \qquad \qquad \rightarrow\text{mass balance substr. conc.}\label{eq:mass_balance_substr}\\
                                          0=& v_{P}'\cdot M_{P}\cdot c_{bio} -c_{P}\cdot v_{bio}'\qquad\,\rightarrow\text{mass balance product conc.}\label{eq:mass_balance_product}\\
                        \boldsymbol{v}' \in& \text{ arg }\underset{\boldsymbol{v}}{\text{max}} 
                            \quad v_{bio}   \qquad \qquad \qquad\,\rightarrow\text{maximize growth}\label{eq:max_growth}\\
                           &  \begin{aligned} \text{s.t.} \quad \boldsymbol{S}\boldsymbol{v}&=\boldsymbol{0} &&\rightarrow\text{mass balances metabolites} \\
                            v_{bio}&\geq f\cdot v_{bio,WT} &&\rightarrow\text{biomass threshold}\\
                           \boldsymbol{v}&\geq \boldsymbol{v^{lower}}\circ (\boldsymbol{B}\boldsymbol{y})&&\rightarrow\text{knockout constraints}\\
                           \boldsymbol{v}&\leq \boldsymbol{v^{upper}}\circ(\boldsymbol{B}\boldsymbol{y}),&&\rightarrow\text{knockout constraints}\label{eq:lower_level_FBA}
                           \end{aligned} 
\end{align}
with $c$ standing for concentration, $v$ being microbial fluxes, and $\boldsymbol{y}\in\{0,1\}^{r}$ denoting the binary knockout vector. 
The units of the variables and parameters can be found in the Nomenclature. 
The size $r$ of this vector corresponds to the number of reactions in the metabolic network. 
In SimulKnock and also throughout this work, we refer to ``knockout'' as ``reaction elimination''.
The subscripts $S$, $P$, $bio$, and $WT$ denote the substrate, the product, the biomass, and the wild-type, respectively.
The biomass flux of the strain before genetic modification $v_{bio,WT}$ is a reference flux and can be either determined experimentally or via FBA.
The superscript apostrophe $'$ formally designates the lower-level variables in the upper level; $i$ denotes the summation index. 
The parameters $\kappa$, $M$, $K_S$, and $f$ stand for the fixed maximum allowable number of knockouts, the molar mass, the Monod constant, and a fixed number between zero and one ensuring a minimal biomass production, respectively. 
The stoichiometric matrix $\boldsymbol{S}$ includes the stoichiometric coefficients of each intracellular reaction, where reversible reactions are split into an irreversible forward and an irreversible backward reaction. 
Accordingly, the mapping matrix $\boldsymbol{B}$ maps reversible reactions with their irreversible counterparts. 

In the upper-level program (Eqs.~\eqref{eq:max_STY}-\eqref{eq:mass_balance_product}), the space-time yield of the fermentation is maximized.
From continuous fermentation and the mass balance of biomass, it follows that the dilution rate equals the growth rate $v_{bio}$, which explains why the growth rate $v'_{bio}$ is part of the space-time yield formula. 
The maximization is subject to two equality constraints, which originate from the mass balances of the substrate (Eq.~\eqref{eq:mass_balance_substr}) and product concentration (Eq.~\eqref{eq:mass_balance_product}).
Further, the upper level contains one inequality constraint (Eq.~\eqref{eq:number_KO}) to ensure that a maximum allowable number of knockouts $\kappa$ (with $\kappa$ as an input parameter) is not exceeded, which is a standard equation from optimal strain design. 
The degrees of freedom are the reaction knockouts $\boldsymbol{y}$ and the substrate feed concentration $c_{S,Feed}$. 

The upper-level program is also subject to the optimality of the lower-level program (Eqs.~\eqref{eq:max_growth}-\eqref{eq:lower_level_FBA}), which is based on FBA and aligns with the lower-level program of OptKnock \citep{Burgard.2003}.
The lower-level objective is to maximize growth subject to the mass balances of the metabolites, a biomass threshold, and the knockout constraints. 
The fluxes $\boldsymbol{v}$ are the degrees of freedom to the lower level. 

\subsection{Bioreactor design for minimized cost}

\begin{figure}[htb]
    \centering
    \begin{subfigure}[b]{0.69\textwidth}
     \centering
        \includegraphics[width=\textwidth]{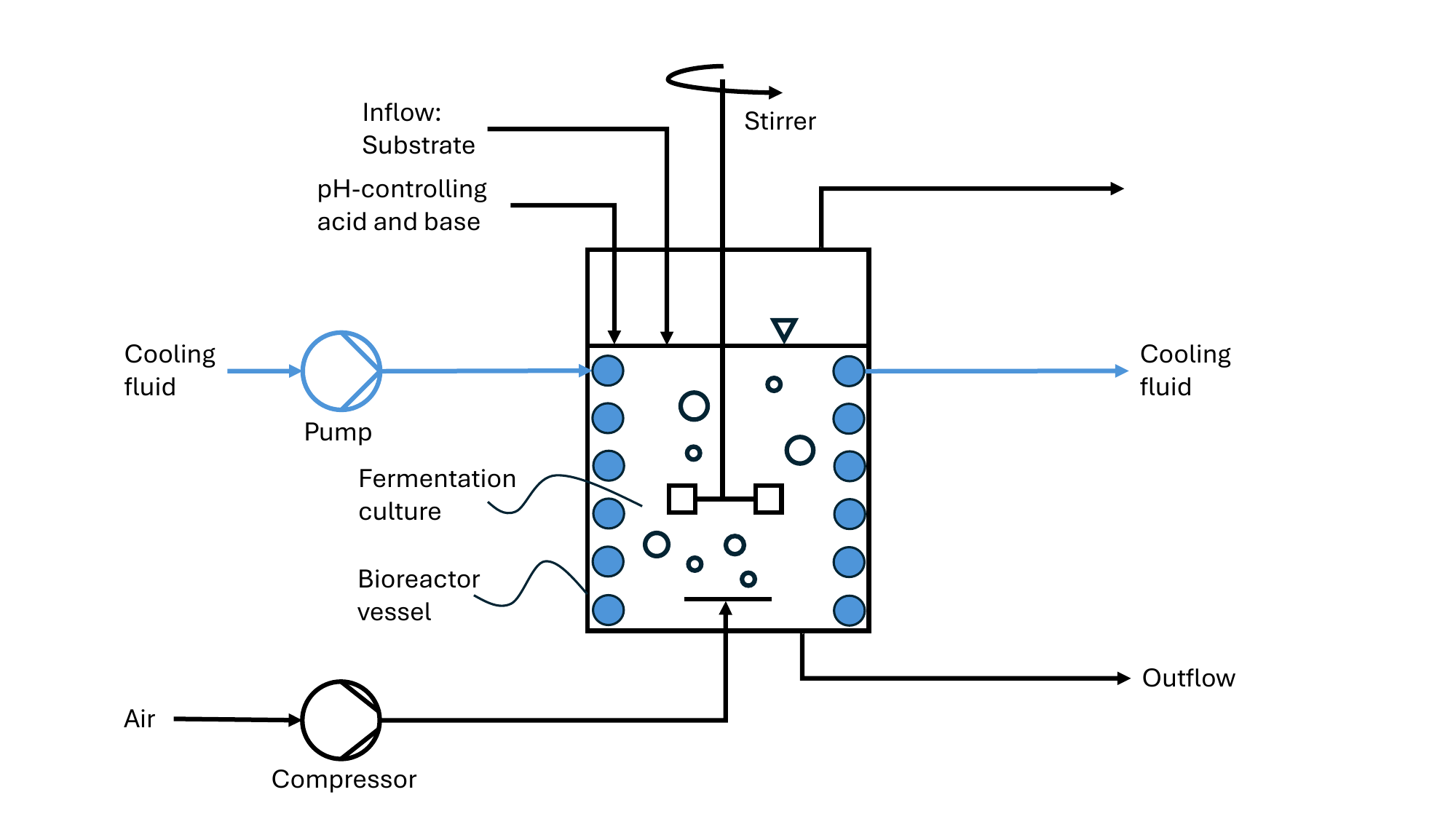}
        \caption{Modeled aspects of a bioreactor}
        \label{subfig:flowsheet_CSTR}
     \end{subfigure}
     \hfill
     \begin{subfigure}[b]{0.29\textwidth}
        \centering
        \includegraphics[trim = {13cm, 0, 6cm, 0}, clip, width=\textwidth]{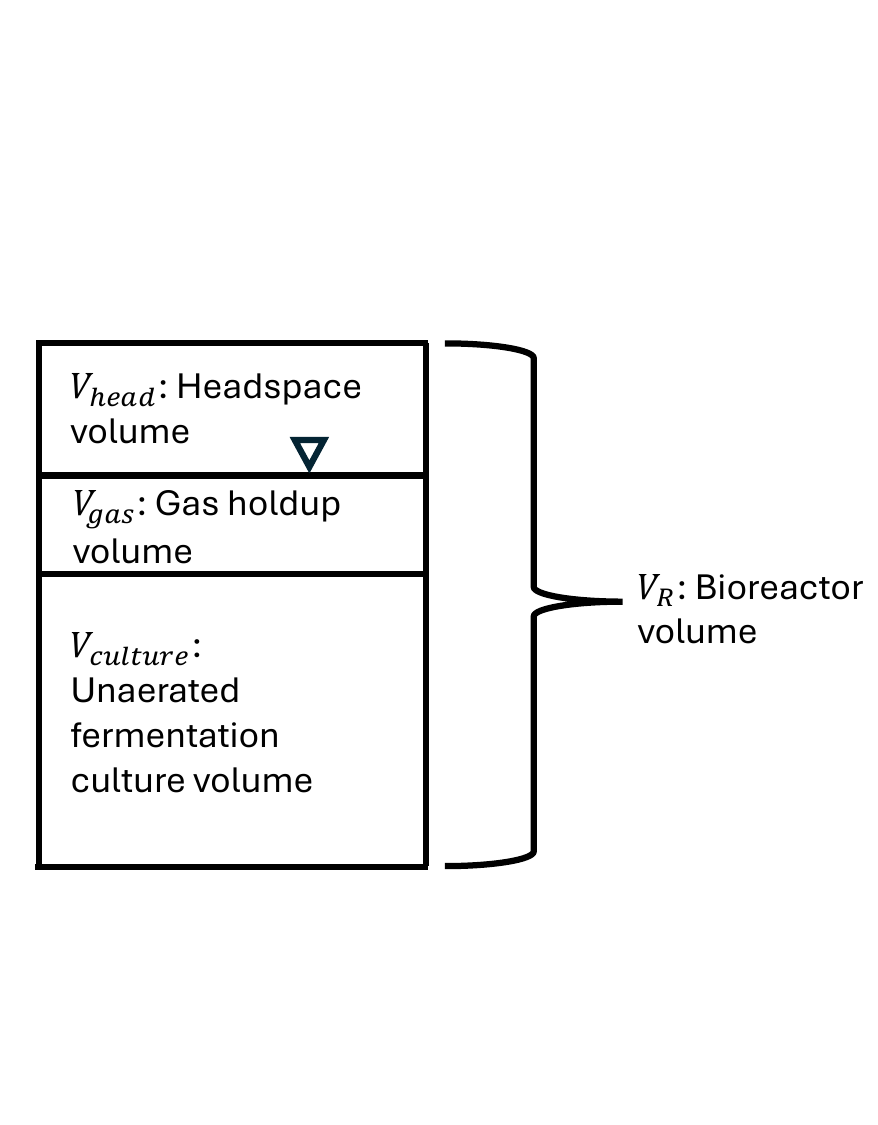}
        \caption{Considered control volumes}
       \label{subfig:bioreactor_volumes}
    \end{subfigure}    
    \caption{Cross section of a continuous stirred-tank bioreactor with the aspects considered in this work in Subfigure a). The cooling system is depicted in blue color. Subfigure b) depicts the different control volumes considered in this work. The unaerated fermentation culture volume and the gas holdup volume were artificially separated in this depiction.}
    \label{fig:Bioreactor_aspects_volumes}
\end{figure}

We modeled a continuous stirred-tank reactor (CSTR) for aerobic fermentation, including pH control, aeration, and cooling. 
The bioreactor model is a set of mass and energy balances, which was combined with cost correlations.
Due to the envisaged connection with the cellular level, cellular variables were introduced into the bioreactor design equations. 
The instances where we performed these replacements will be highlighted in the text. 
Figure~\ref{fig:Bioreactor_aspects_volumes} depicts the aspects that were considered for bioreactor design.

We focused on the bioreactor, i.e., neither considered upstream nor downstream operations. 
Upstream operations would include sterilization of the substrate, media preparation, or filtration of the air. 
Downstream operations would include the purification of the product. 
Considering the bioreactor represents a first step, as adding upstream and downstream units would further increase the complexity.
Furthermore, the fermentation may represent the biggest OPEX cost factor in large-scale operation \citep{Nieto.2020}. 
Moreover, we did not consider operating labor costs. 
In Figure~\ref{subfig:bioreactor_volumes}, the considered control volumes within the bioreactor are depicted. 
We ignored the volume of the baffles, the mixer, the sparger, and other internal tubing because they are small compared to the presented control volumes in the production scale.  

The economic viability of a bioprocess depends on its production cost.
We minimize the hourly production cost, which is composed of CAPEX and OPEX. 
We assumed CAPEX to be represented by the investment cost of the bioreactor and the compressor; while we assumed the capital cost of the pump to be negligible.
For calculating OPEX, we considered raw material costs of the cooling agent, the substrate, and the pH control; as well as power.
We considered air for aeration and assumed that air is free of charge.
Further, we assumed that due to the continuous operation, no sterilization and maintenance of the reactor is needed.

\begin{align}
    \text{Cost} =& \underbrace{\frac{C_R + C_{compress}}{t_{annual}\cdot t_{amortization}}}_\text{Investment} + \underbrace{\frac{\dot{m}_{cooling}}{\rho_{cooling}} C_{cooling}}_\text{Cooling agent} \nonumber\\
    &+ \dot{V} (\underbrace{C_S\,c_{S,Feed}}_\text{Substrate}+\underbrace{C_{acid}\,M_{acid}\,c_{acid} + C_{base}\,M_{base}\,c_{base}}_\text{pH control} ) \nonumber\\
    &+\underbrace{C_{Power}\,\left(P_{compress} + P_{agitation} + P_{cooling}\right)}_\text{Power}, 
\label{eq:cost}
\end{align}
where $t$ denotes the time, $C$ denotes the specific costs; $\dot{m}$ stands for mass flow; $\rho$ stands for density,  and $\dot{V}$ stands for the volumetric hourly reactor outflow.
The power input is denoted with $P$.
The subscript $R$ stands for the reactor, and $compress$ stands for compression, whereas $acid$ and $base$ in the subscript denote the acid and base of the pH control, respectively.
Sulfuric acid was chosen as the acidic pH control component, and sodium hydroxide as the basic pH control component.
The specific costs $C_S$, $C_{acid}$, $C_{base}$, and $C_{Power}$ as well as $t_{amortization}$ and $t_{annual}$ are parameters; the values that were used for this work are tabulated in Table~\ref{tab:bioreactor_parameter_values}. 

The reactor investment cost was calculated based on the principle of economy of scale \citep{Powell.2008, Perry.2008} in Eq.~\eqref{eq:invest_reactor}.
Therein, we ignored the inflation factor since we assumed that inflation is balanced out by cheaper production. 
Further, we ignored CAPEX for media tanks, filters, or seed tanks, which can be associated with the reactor \citep{Negulescu.2023}.
The investment cost for the compressor was calculated using the equation from \citet{Luyben.2018} in Eq.~\eqref{eq:invest_compressor}, based on Guthrie's correlations \citep{Douglas.1988}.
\begin{align}
    C_{R} =& C_{R,ref} \left(\frac{V_R}{V_{R,ref}}\right)^{0.6} \label{eq:invest_reactor}\\
    C_{compress} =& 5840\, P_{compress}^{0.82}\,, \label{eq:invest_compressor}
\end{align}
 where the subscript $ref$ denotes a reference state.
 Below the reference volume $V_{R,ref}$, the hourly cost was set equal to $C_{R,ref}$.
Again, the chosen values for the parameters, that is, $C_{R,ref}$ and $V_{R,ref}$, can be found in Table~\ref{tab:bioreactor_parameter_values}.

In Eq.~\eqref{eq:reactor_volume_to_broth_volume}, the reactor volume was related to the fermentation culture volume. 
The hourly reactor outflow $\dot{V}$ was related to the culture volume and was calculated from the annual production capacity \citep{GarciaOchoa.2011} by:
\begin{align}
    V_R &= \beta \cdot (V_{culture} + V_{gas}) \label{eq:reactor_volume_to_broth_volume} \\
    V_{culture} &= \frac{\dot{V}}{v_{bio}}\label{eq:reactor_outflow_to_volume}\\
    \dot{V} &= \frac{m_{P,annual}}{t_{annual} \cdot c_P} \label{eq:reactor_outflow_to_production_capacity}\\
    V_{culture} &= \frac{\pi}{4} D^2 H_{culture} \label{eq:culture_volume}\\
    H_{culture} &= q \cdot D,\label{eq:culture_height}
\end{align}
where $V$ denotes volume and the produced mass is denoted with $m$; $D$ for diameter, and $H$ for the height of the fermentation culture.
The multiplier $\beta$ accounts for the reactor volume $V_R$ being larger than the culture volume $V_{culture}$ due to a headspace volume \citep{VantRiet.1991} and the multiplier $q$ relates the height to the diameter of the culture (see Table~\ref{tab:bioreactor_parameter_values} for chosen values). 
In Eq.~\eqref{eq:reactor_outflow_to_volume}, the dilution rate was replaced by the growth rate $v_{bio}$, due to continuous fermentation conditions.
With the growth rate being a variable of the cellular level, a bioreactor parameter was replaced with a cellular variable.

As comes clear from Eqs.~\eqref{eq:culture_volume}-\eqref{eq:culture_height}, we assumed the fermentation culture to be of cylindrical shape, just like the reactor.
Together with Eq.~\eqref{eq:reactor_volume_to_broth_volume}, the size of the reactor is determined during optimization, which is referred to as optimal bioreactor design.

\subsubsection*{Numbering-up of reactors}
To take into account large production capacities with a maximum standard tank size \citep{Powell.2008} and to avoid large spatial distributions (dead zones), we implemented an optional parallel setup of reactors, called numbering-up.
For numbering-up, Eq.~\eqref{eq:numbering_up} instead of Eq.~\eqref{eq:reactor_outflow_to_volume} is chosen. 
\begin{equation}
    V_{culture} \cdot n_R = \frac{\dot{V}}{v_{bio}}, \label{eq:numbering_up}
\end{equation}
where $n_R$ is the number of reactors set up in parallel. 
An upper bound on the number of reactors was set to $n_R^{upper} = 10$; an upper bound on the reactor size was set to $V_R^{upper} = \SI{10}{\cubic\meter}$.
Only at a production capacity of \SI{e6}{\kg\per\year}, the maximum number of reactors was set to 100 and the maximum reactor volume to \SI{100}{\cubic\meter}, respectively. 
Whether the numbering-up option becomes active depends on the production rate chosen. 

In the following, the terms for pH control, aeration, and cooling are described in detail.

\subsubsection*{pH control}
The pH level is a decisive factor for the growth and productivity of a microbial strain \citep{Parhad.1974} and, hence, pH control is an integral part of a bioprocess. 
The pH in fermentation cultures can be distinguished into extracellular and intracellular pH \citep{Straathof.2023}. 
In this work, we focused on the extracellular pH.
The extracellular pH impacts different aspects of a bioprocess, e.g., growth, product formation, and maintenance requirements in the cell \citep{Straathof.2023}. 
Of these impacts, we chose to model the impact of pH on growth for three reasons. 
First, even if exemplary publications exist, e.g., for modeling maintenance requirements \citep{Du.2019}, strain-wide data are not available yet. 
Finding a general model for the influence of pH which applies to several strains and several products proved to be difficult. 
Second, the optimization would be further restricted, for example, by implementing a Luedeking-Piret equation or similar for inhibited product formation. 
Third, like SimulKnock, this work is also based on growth-coupled production, which already links growth with product formation via the GEM. 
By varying the pH and modeling the impact on growth, we assumed that the predicted genetic modifications would not affect the organism's behavior toward changing pH conditions. 

In the absence of experimental data, the extracellular pH was modeled based on the simplified definition of pH (Eq.~\eqref{eq:def_pH}) and on the idea of a charge balance by \citet{Campos.2003}.
We conducted preliminary studies, which were in alignment with \citet{Silverstein.2000}, who argued that the impact of weak ions is usually smaller than the impact of strong ions. 
Hence, we only considered the strong ion balance for pH calculation as follows
\begin{align}
    pH =& -\log_{10}\,c_{H^+} \label{eq:def_pH}\\
    c_{H^+} =& \underbrace{\frac{c_{bio}}{v_{bio}}\left(\sum_{i}^{n_{acid, strong}}v_i - \sum_{i}^{n_{base, strong}}v_i\right)}_\text{strong ion balance}  + \underbrace{c_{acid} - c_{base}}_\text{pH control}\, , \label{eq:strong_ion_balance}
\end{align}
where $c_{H^+}$ denotes the equilibrium concentration of hydrogen ions.

Eq.~\eqref{eq:strong_ion_balance} shows that we replaced the concentrations of the acid and base in the fermentation culture with their respective fluxes from the cellular level such that $c_i = \frac{c_{bio}}{v_{bio}} \cdot v_i$. 
As the extracellular pH was modeled, only fluxes of so-called exchange reactions were taken into account. 
Exchange reactions denote reactions in the GEM that connect the extracellular with the intracellular space.

The impact of pH on growth can be modeled by adapting kinetic models \citep{Nath.2011, Straathof.2023, Wang.2009}; however, they mostly require experimental data.  
We introduced an inhibition factor \citep{Zwietering.1993}, depicted in Eq.~\eqref{eq:impact_pH_growth} using the cardinal pH model (CPM) \citep{Rosso.1995} in Eq.~\eqref{eq:CPM_model}. 

\begin{align}
    v_{bio}^{max} =& v_{bio}^{opt} \gamma\qquad \text{ for all } pH^{min} < pH < pH^{max}\label{eq:impact_pH_growth}\\
    \gamma =& \frac{(pH - pH^{min}) (pH - pH^{max})}{(pH - pH^{min}) (pH - pH^{max}) - (pH - pH^{opt})^2}, \label{eq:CPM_model}
\end{align}
where $\gamma$ is the growth inhibition factor, and $pH^{min}$, $pH^{max}$, and $pH^{opt}$ indicate strain-specific parameters for the minimum, the maximum, and the optimum pH, respectively. 
The growth rate $v_{bio}^{opt}$ denotes the maximum growth rate at optimal pH, which corresponds to the parametric value of $v_{bio}^{max}$ in SimulKnock (Eq.\eqref{eq:mass_balance_substr}).
Below $pH^{min}$ and above $pH^{max}$, the growth rate $v_{bio}$ is zero. 
Since growth rates equal to or close to zero are techno-economically not sensible, we set upper and lower bounds on the pH variable, namely $pH^{lower}$ and $pH^{upper}$. 
These bounds are tighter than the parameters $pH^{min}$ and $pH^{max}$, respectively, and represent the vital zone of the organism.
As a side-effect, by using these bounds, we could avoid implementing a logical constraint resulting from Eq.~\eqref{eq:impact_pH_growth}. 
We chose the CPM model due to its wide applicability to many microbial strains and due to its small number of additional parameters, i.e., $pH^{min}$, $pH^{max}$, and $pH^{opt}$, which can be found in the literature. 

To sum up, we modeled the pH based on the exchange fluxes of the cellular level, with a focus on the strong ions. 
The cost for pH control was taken into account by the required concentration of pH-regulating acid and base. 
Further, we modeled the impact of pH on the growth by introducing an inhibition factor to the Monod kinetics.
In the next step, the oxygen supply is examined in detail. 

\subsubsection*{Oxygen supply via aeration and agitation}

In aerobic fermentation, oxygen is a key component to the fermentation that may control the metabolism \citep{GarciaOchoa.2010}.
Due to its low solubility in water, oxygen has to be supplied continuously to the fermentation culture and its uniform distribution is a challenge of bioreactor design and scale-up \citep{GarciaOchoa.2009, Clarke.2013, Liu.2020}.
In this work, we assumed an ideal mixture of the fermentation culture. 
The supply of oxygen was modeled via compressed air through a sparger and the distribution of oxygen via stirring \citep{PortodeSouzaVandenberghe.2022}. 
For modeling the required electrical power, we considered the compression of air for aeration and the mixing for agitation.

The mass balance for dissolved oxygen in the fermentation culture reads 
\begin{equation}
    \frac{\mathrm{d}c_{O_2}}{\mathrm{d}t} = OTR - OUR = 0, \label{eq:mass_balance_oxygen}
\end{equation}
where $OTR$ is the Oxygen Transfer Rate and $OUR$ is the Oxygen Uptake Rate. 
We assumed a pseudo-steady state of the oxygen transfer, as it is practiced by other modeling examples, for example, \citet{Humbird.2017}. 
This assumption results in $\frac{\mathrm{d}c_{O_2}}{\mathrm{d}t} = 0$, from which follows $OUR = \,OTR$.

The OUR was calculated directly by its correlation with the oxygen uptake flux $v_{O_2}$ in Eq.~\eqref{eq:OUR}, based on the definition of OUR \citep{Liu.2020,GarciaOchoa.2011}.
By directly using the oxygen uptake flux, we replaced a parameter of the bioreactor level, namely, the respiration rate, with a variable from the cellular level, which is $v_{O_2}$. 
Thus, the predictions at the bioreactor level profit from the accuracy of the GEM. 
In Eq.~\eqref{eq:oxygen_gas_flow}, we assumed that the provided oxygen gas stream equals the volume flow of oxygen taken up. 
We assumed air to be an ideal gas, with a volumetric composition of 21\% oxygen (see Eq.~\eqref{eq:air_gas_flow}).
\begin{align}
    OUR =& \, v_{O_2} M_{O_2} c_{bio}\label{eq:OUR}\\
    \dot{V}_{O_2} =& \frac{OUR \cdot V_{culture}}{\rho_{O_2}} \label{eq:oxygen_gas_flow}\\
    \dot{V}_{air} =& \frac{\dot{V}_{O_2}}{0.21}. \label{eq:air_gas_flow}
\end{align}

As it is common practice, based on the two-film theory, the resistance in the liquid boundary layer of the gas-liquid interface was assumed to be decisive for the OTR \citep{Clarke.2013, Seidel.2021}. 
Hence, the resistances in the gaseous boundary layer and the boundary layer at the cell membrane were omitted (Eq.~\eqref{eq:OTR_kLa}). 
The $k_La$-value was then calculated using the van't Riet correlation \citep{VantRiet.1979} for a non-coalescing medium in Eq.~\eqref{eq:vantRiet}.
The power of agitation is taken into consideration through this equation.
Further, we ignored the influence of aeration on power input because this depends on specific bioreactor design decisions, e.g., type and amount of agitators. 
In Eq.~\eqref{eq:superficial_gas_velocity}, the superficial gas velocity was calculated based on the assumption of a cylindrical reactor (Eq.~\eqref{eq:surface_cross_section}).
\begin{align}
    OTR =& k_La\left( c_{O_2,equil} - c_{O_2}\right) = k_La \cdot (1- \varphi) \cdot c_{O_2,equil} \label{eq:OTR_kLa}\\
    k_La =& 0.002 \left( \frac{P_{agitation}}{V_{culture}} \right)^{0.7} u_s^{0.2} \label{eq:vantRiet} \\
    u_s = & \frac{\dot{V}_{air}}{A_{cross}}\label{eq:superficial_gas_velocity} \\
    A_{cross}  =& \frac{\pi}{4} D^2, \label{eq:surface_cross_section}
\end{align}
where $k_La$ is a lumped parameter for overall volumetric oxygen transfer. 
The concentration $c_{O_2}$ refers to the oxygen concentration in the fermentation culture (also called dissolved oxygen), whereas $c_{O_2,equil}$ denotes the oxygen concentration at the gas-liquid interface, i.e., the solubility of oxygen in the fermentation culture when assuming chemical equilibrium at the interface. 
The parameter $\varphi$ correlates oxygen concentration in the culture with oxygen solubility. 
In Eq.~\eqref{eq:OTR_kLa}, we set $c_{O_2} = \varphi \cdot c_{O_2,equil}$ with $\varphi = 0.3$ \citep{Bernard.2001,Strandberg.1991} to avoid oxygen limitation of the fermentation, which can occur already at low dissolved oxygen levels \citep{Hansen.2022}.   
Further, $u_S$ denotes the superficial gas velocity.
The cross-sectional area of the fermentation culture $A_{cross}$ equals the base surface of the fermentation culture.

With the superficial gas velocity and the agitation power at hand, the gas holdup can be calculated. 
The gas holdup specifies the dispersed gas in an aerated bioreactor.
It becomes more and more important with increasing bioreactor sizes and commonly lies between 10 and 30\% \citep{Villadsen.2011, Nienow.2009}. 
The gas holdup can be correlated to the superficial gas velocity and the specific agitation power input \citep{Sieblist.2010,Yawalkar.2002, Bailey.1986} and was calculated in Eq.~\eqref{eq:correlation_gas_holdup}, based on the correlation of \citet{Weissgarber1976}.
The gas holdup volume was calculated using the definition of the gas holdup in Eq.~\eqref{eq:def_gas_holdup}.
\begin{align}
    \epsilon_{gas} =& 1.12 \cdot \left( \frac{P_{agitation}}{V_{culture}\cdot \rho_{culture}}\right)^{0.29} \cdot u_{S}^{0.6}\label{eq:correlation_gas_holdup}\\
    \epsilon_{gas} =& \frac{V_{gas}}{V_{gas}+V_{culture}}, \label{eq:def_gas_holdup}
\end{align}
where $\epsilon_{gas}$ stands for the gas holdup. 
The culture density was approximated with the density of water \citep{Mayer.2023, MartinezSalas.1981}.

The required power for aeration was calculated by the thermodynamic equation of adiabatic, isentropic compression of ideal gases \citep{Biegler.1997, Ulonska.2016, Cardoso.2020} in Eq.~\eqref{eq:power_compressor}.
We assumed the feedstock cost of air to be small compared to compressor cost and cost for power, which is why we neglected them.
\begin{align}
P_{compress} =& \dot{n}_{air} \frac{1}{\eta_{compress}} \frac{\alpha}{\alpha -1} R T \left( \left( \frac{p_{air}}{p^0} \right ) ^{\frac{\alpha-1}{\alpha}} -1 \right) \nonumber\\
=& p^0 \dot{V}_{air} \frac{1}{\eta_{compress}} \frac{\alpha}{\alpha -1} \left( \left( \frac{p_{air}}{p^0} \right ) ^{\frac{\alpha-1}{\alpha}} -1 \right) \label{eq:power_compressor},
\end{align}
where $\dot{n}$ denotes the mole flow rate; $\eta$ is an efficiency number, $\alpha$ is a polytropic exponent.
The ideal gas constant $R$ and the temperature $T$ cancel each other out when applying the ideal gas law $p^0 \dot{V}_{air}= \dot{n}_{air} R T$.
Atmospheric pressure is denoted with $p^0$.
The chosen values for the parameters $\eta_{compress}$, $\alpha$, and $p_{air}$ are tabulated in Table~\ref{tab:bioreactor_parameter_values}.
Note that $p_{air}$ was set to \SI{2.5}{\bar}.

In summary, aeration was modeled by setting the oxygen uptake at the cellular level about the required oxygen supply and the oxygen transfer rate. 
The oxygen transfer rate is decisive for the required agitation power. 
The cost of aeration was calculated by taking into account the capital cost of the compressor and the required power for compression. 
The next and last modeling aspect is an energy balance including cooling.

\subsubsection*{Cooling}
In aerobic fermentation, cooling is an important aspect of bioreactor design, since temperature affects the microbial metabolism and the oxygen solubility \citep{Clarke.2013}.
Especially for large bioreactors, providing enough heat transfer surface is crucial, since the ratio of surface area to reactor volume decreases with increasing reactor size \citep{Clarke.2013}. 
Heat is mainly transferred from the metabolic activity and the agitation \citep{Rao.2006}, as depicted in Eq.~\eqref{eq:total_heat_generation}, with the agitation heat equalling the agitation power input \citep{Knoll.2005} in Eq.~\eqref{eq:heat_agitation}.
Again, we ignored the influence of aeration on power input. 
\begin{align}
    \dot{Q} =& \dot{Q}_{agitation} + \dot{Q}_{met} \label{eq:total_heat_generation}\\
    \dot{Q}_{agitation} =& P_{agitation} ,\label{eq:heat_agitation}
\end{align}
where $\dot{Q}$ denotes the heat flow; and the subscript $met$ stands for metabolic.

In bioreactor design, it is common to correlate metabolic heat with the OUR or the OTR \citep{Doran.1995,Knoll.2005,Cooney.1969}.
We applied the correlation 
\begin{equation}
    \dot{Q}_{met} = \delta \cdot \frac{OUR \cdot V_{culture}}{ M_{O_2}}, 
    \label{eq:correlation_metabolic_heat}
\end{equation}
where the correlation parameter  $\delta = \SI{0.13}{\kilo\watt\hour\per\mol}$ is a mean value of the values specified in the aforementioned publications. 

We assumed that the temperature of the fermentation culture is known, constant in time, and uniform in space. 
The heat transfer through the cooling agent was modeled in Eq.~\eqref{eq:heat_transfer_cooling_fluid}.
The required power for pumping the cooling agent through the cooling coil was calculated similarly to \citet{Knoll.2005} in Eq.~\eqref{eq:power_cooling}.
\begin{align}
    \dot {Q} =& \, \dot{m}_{cooling} \cdot c_{p,heat} (T_{cooling}^{in} - T_{cooling}^{out}) \label{eq:heat_transfer_cooling_fluid}\\
    P_{cooling} =& \frac{\dot{m}_{cooling} \cdot g \cdot (H_{culture}+H_{gas})}{\eta_{pump}}, \label{eq:power_cooling}
\end{align}
where $c_{p,heat}$ is the specific heat capacity of the cooling agent at the mean temperature of the cooling agent $T_{mean}$ and $g$ denotes gravity. 
We chose cooling water as the cooling agent at input temperature $T_{cooling}^{in} = \SI{15}{\celsius}$.
The final temperature of the cooling water $T_{cooling}^{out}$ was assumed to be \SI{25}{\celsius}, which was \SI{5}{\celsius} below the temperature of the fermentation culture $T_{culture}$ \citep{Heinzle.2006}. 
The mean temperature was calculated by the logarithmic mean temperature difference, i.e., $T_{mean} = \SI{20.9}{\celsius}$.

\subsection[SimulKnockReactor]{SimulKnockReactor: Combined optimization formulation}

We combined the formulation of SimulKnock with the equations from bioreactor design and cost analysis. 
The combined optimization formulation SimulKnockReactor reads: 
\begin{align}
    \textbf{min } &\textbf{cost} \nonumber\\
    \text{s.t. } & \text{cost calculation (Eq.~\eqref{eq:cost}),} \nonumber \\
                & \text{investment cost reactor and compressor (Eqs.~\eqref{eq:invest_reactor}-\eqref{eq:invest_compressor}),}\nonumber\\
                & \text{bioreactor sizing (Eqs.~\eqref{eq:reactor_volume_to_broth_volume}-\eqref{eq:culture_height} and \eqref{eq:correlation_gas_holdup}-\eqref{eq:def_gas_holdup})}, \nonumber\\
                & \text{pH control (Eqs.~\eqref{eq:def_pH}-\eqref{eq:strong_ion_balance}}), \nonumber \\
                & \text{growth inhibition through pH (Eqs.~\eqref{eq:impact_pH_growth}-\eqref{eq:CPM_model}),} \nonumber \\
                & \text{oxygen uptake and transfer (Eqs.~\eqref{eq:mass_balance_oxygen}-\eqref{eq:OTR_kLa}),} \nonumber\\
                & \text{power for agitation \& compression (Eqs.~\eqref{eq:vantRiet}-\eqref{eq:surface_cross_section} and \eqref{eq:power_compressor}),} \nonumber \\
                & \text{energy balance and heat transfer (Eqs.~\eqref{eq:total_heat_generation}-\eqref{eq:correlation_metabolic_heat} and \eqref{eq:heat_transfer_cooling_fluid}),} \nonumber \\
                &\text{pumping power for cooling (Eq.~\eqref{eq:power_cooling}),} \nonumber \\
                & \text{mass balances fermentation \& number of knockouts (Eqs.~\eqref{eq:number_KO}-\eqref{eq:mass_balance_product}),} \nonumber  \\
                &\textbf{max growth rate} \nonumber \\
                &\text{s.t. mass balances \& knockout constraints (Eqs.~\eqref{eq:lower_level_FBA}).} \nonumber
\end{align}

The difference to SimulKnock lies in the upper-level program: it includes bioreactor design aspects now. 
Thus, the objective function is different. 
The number of modeled compounds and the degrees of freedom have increased. 
A detailed comparison between SimulKnock and this work's extension, SimulKnockReactor, is depicted in Table~\ref{tab:comparison_SimulKnock_Extension}.

\begin{table}[h!tb]
    \centering
    \caption{Comparison between the formulations SimulKnock and the extension SimulKnockReactor of this work. UL: upper level, LL: lower level, FBA: flux balance analysis, conc.: concentration, KO's: knockouts}
    \label{tab:comparison_SimulKnock_Extension}
     \resizebox{\textwidth}{!}{
    \begin{tabular}{lll}
    \hline
     Aspect    & SimulKnock& This work: SimulKnockReactor    \\
    \hline
    Optimization class & Bilevel (UL \& LL) & Bilevel (UL \& LL) \\
    Input parameters & Metabolic network, & Metabolic network, \\
    &maximum number of KO's & maximum number of KO's,\\
    && production capacity \\
    Implemented kinetics & Monod or  & Monod with pH-inhibition \\
            & Michaelis-Menten &  \\
    \hline
    Upper Level (UL) & Fermentation & Fermentation \& bioreactor  \\
    Objective function & Space-time yield & Bioreactor cost \\
    Modeled compounds & Substrate, product,  & Substrate, product, biomass, \\
            &biomass&oxygen, hydrogen, strong ions \\
    Degrees of freedom &  Knockout variable, & Knockout variable,\\
                            & substrate feed conc. & substrate feed conc., \\
                            && conc. pH control agents,\\
                            && agitation power, \\
                            && compression power, \\
                            && mass flow cooling agent,\\
                            && optional: number of parallel reactors \\
    Additional variables & Growth/ dilution rate & Growth/ dilution rate,\\
                &&height and diameter of \\
                &&fermentation culture, \\
                && volume of reactor,\\
    
    &  & compression of air, \\
            &&agitation, \\ 
            &&cooling, \\
            &&pH control \\
    \hline
    Lower Level (LL) & FBA &FBA \\
    Objective function & Growth rate & Growth rate \\ 
    Degrees of freedom & Metabolic fluxes & Metabolic fluxes \\
    \hline
    \end{tabular}
    }
\end{table}

\subsection{Implementation}
For implementation and solution, we reformulated the bilevel program to a single-level program using strong duality. 
This is a standard technique also used in OptKnock \citep{Burgard.2003} and in SimulKnock \citep{Ziegler.2024}, which comes at the expense of an increased number of variables due to the introduction of so-called dual variables.
The resulting optimization program is a mixed-integer nonlinear optimization program. 
Similar to SimulKnock, we used auxiliary variables to reformulate the nonlinear constraints of the substrate mass balance with Monod kinetics (Eq.~\eqref{eq:mass_balance_substr}) and the CPM model (\eqref{eq:CPM_model}) such that they could be solved by commercial solvers, e.g., Gurobi \citep{GurobiOptimization.2023}.
The remaining nonlinear functions, i.e., the investment cost calculation of the reactor (Eq.~\eqref{eq:invest_reactor}) and the compressor (Eq.~\eqref{eq:invest_compressor}), the definition of pH (Eq.~\eqref{eq:def_pH}), the van't Riet correlation of the $k_La$-value (Eq.~\eqref{eq:vantRiet}), and the gas holdup correlation (Eq.~\eqref{eq:correlation_gas_holdup})  are treated through dynamic outer approximations as part of the branch and bound tree. 
To achieve this, we used so-called general constraint functions from Gurobi \citep{GurobiOptimization.2023}. 
The reformulated optimization formulation was implemented with the software Pyomo \citep{hart2011pyomo,bynum2021pyomo} and solved with Gurobi v11.0.3 \citep{GurobiOptimization.2023}, using the preprocessing procedure of SimulKnock \citep{Ziegler.2024}.

\section{Case Studies}
In the following, we show the applicability and the mode of action of SimulKnockReactor in exemplary case studies. 
We applied SimulKnockReactor for the production of three different target chemicals:  formate, acetate, and succinate. 
We chose the well-examined \textit{E. coli} K12 as an exemplary organism and embedded two different metabolic networks of the organism: the \textit{E. coli} core network \citep{Orth.2010a} and the GEM \textit{i}ML1515 \citep{Monk.2017}. 
The  optimization problem with the embedded GEM (single-level reformulation) is numerically very challenging, due to the large number of variables (around \SI{1.3e4}{} continuous variables, \SI{2.7e3}{} binary variables, and one integer variable), large number of constraints (\SI{1.5e4}{}), and the nonlinearity (7 nonlinear equality constraints).
For all optimization runs, we limited the runtime to \SI{10}{\hour} on 48 cores. 
When the optimization does not converge to a global solution within this runtime, we read out the local solution.
We only reported available solutions and further discriminated against solutions having a higher objective value than the core network solution.
The solution with the embedded core network is supposed to be worse since the core network has fewer reaction pathways included.
In general, the local solution is an upper bound of the global solution; hence, better solutions might exist. 

In the following, we examined the influence of the number of knockouts and the production capacity on the total cost and the individual cost factors.
Moreover, we compared our simultaneous approach with a sequential approach and to a thermochemical formate production process from industry. 

 \subsection[Cost over number of knockouts] {Decreasing cost with an increasing number of knockouts}
 \label{sec:cost_over_knockouts}
\begin{figure}[h!tb]
     \centering
     \begin{subfigure}[b]{0.49\textwidth}
        \centering
        \includegraphics[width=\textwidth]{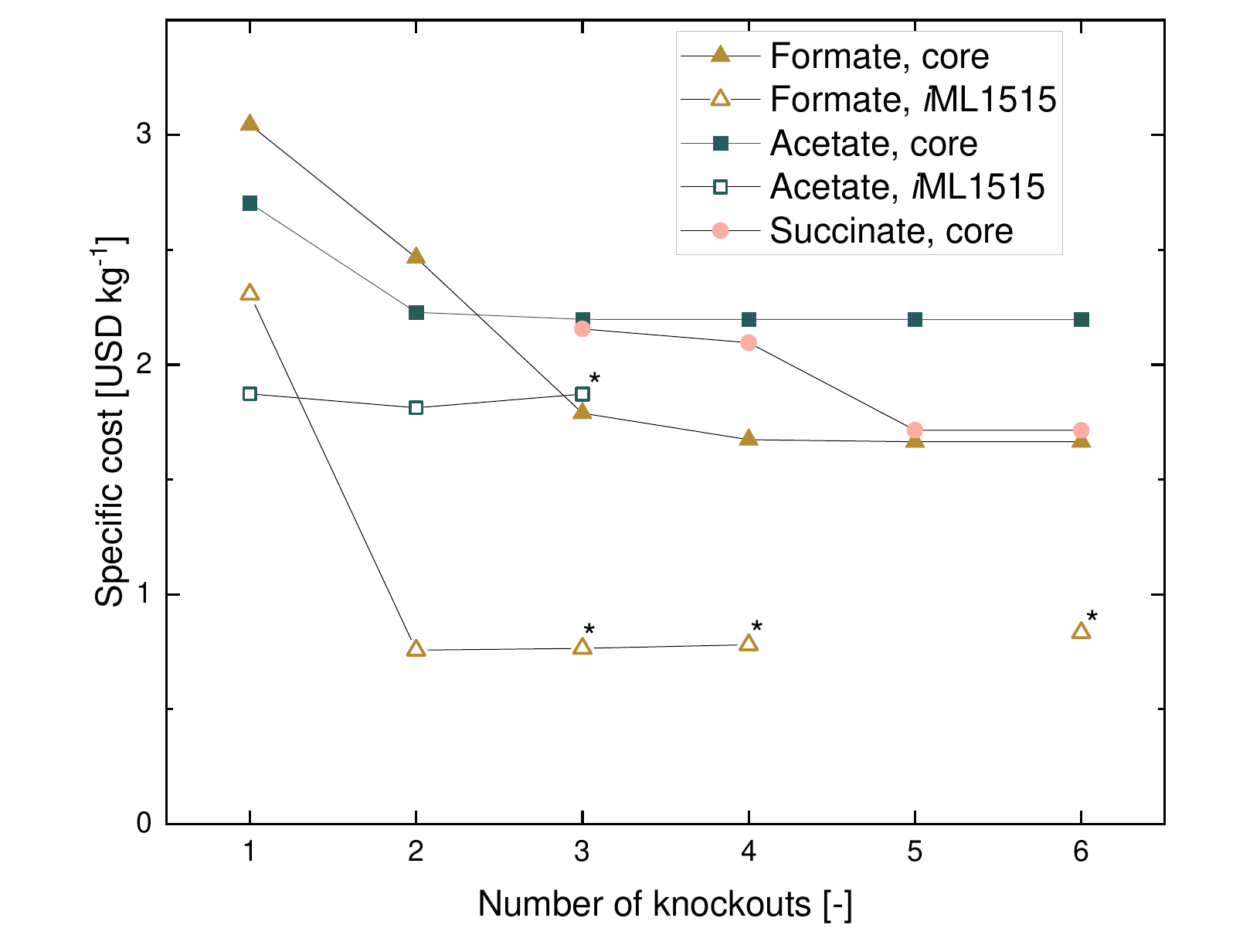}
        \caption{Total cost}
        \label{subfig:cost_over_knockouts_total}
    \end{subfigure}
    \hfill
     \begin{subfigure}[b]{0.49\textwidth}
     \centering
         \includegraphics[width=\textwidth]{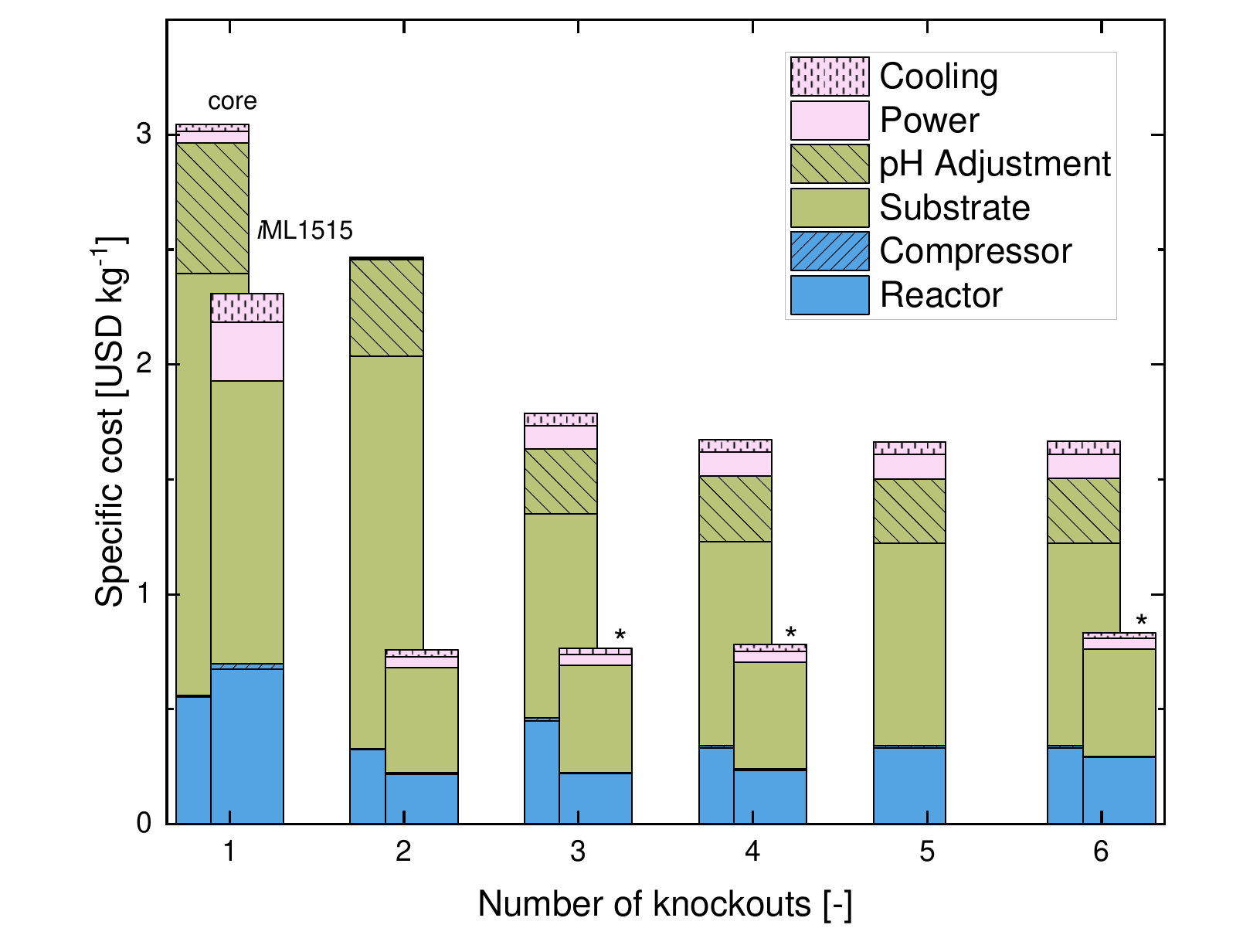}
         \caption{Formate: cost factors}
         \label{subfig:cost_over_knockouts_formate}
     \end{subfigure}
     \caption{Influence of the maximum allowable number of knockouts on the specific production cost. The production capacity was set to \SI{e3}{\kg\per\year}. The asterisk * indicates locally optimal solutions; all other solution points are globally optimal. The embedded networks are \textit{E. coli} core \citep{Orth.2010a} (indicated or bars in the back) and \textit{i}ML1515 \citep{Monk.2017} (indicated or bars in the front) USD: US-Dollar. (a) Total cost for acetate, succinate, and formate, (b) Formate production cost divided into the individual cost factors. Blue, green, and pink color denotes investment cost, raw materials, and utilities, respectively. }
     \label{fig:cost_increasing_KOs}
 \end{figure}

In this case study, we increased the total number of maximum allowable number of knockouts from one to six. 
The total specific cost for formate, acetate, and succinate are depicted in Figure~\ref{subfig:cost_over_knockouts_total}.

The total cost decreases with an increasing number of knockouts. 
The slight increase in cost for acetate and formate production with the embedded GEM can be explained by the fact that solutions responsible for the rise are only local optima. 
A lower (globally optimal) solution value most likely exists. 
For succinate production, the optimization software, that is, the solver, only returns a feasible solution for three knockouts or more. 
For succinate production with the embedded GEM, the solver did not return any global or local solution. 

For each target product, there is a plateau at which a further increase in the number of knockouts does not reduce costs any further or only insignificantly.
For example, this plateau is reached for succinate production at five knockouts and for formate production with the embedded GEM at two knockouts. 

Where a solution is available, the costs with the embedded GEM are lower than the costs with the embedded core network.
For one maximum allowable knockout, the predicted knockout for acetate production with the embedded core network and the GEM is similar, namely ATP synthase. 
The same holds for formate production. 
Presumably, the GEM includes a more efficient pathway toward biomass, which can explain the better objective values with the embedded GEM despite the same knockout strategy. 
For two maximum allowable knockouts, the predicted knockouts with the embedded GEM and core network differ. 
In formate production with the embedded GEM, for example, the knockout of ATP synthase and glutamate dehydrogenase results in a drastic cost decrease. 
For comparison, predicted knockouts with the embedded core network are phosphoenolpyruvate carboxylase and NADH dehydrogenase. 
Hence, differing knockout predictions are a second explanation for different objective values. 
Moreover, with the embedded GEM, the plateau is reached earlier. 
For example, acetate costs stagnate after one knockout with the embedded GEM, while they reach a plateau only after two knockouts with the embedded core network. 
This could also be due to the decreased pathway options in the core network. 

Exemplarily, Figure~\ref{subfig:cost_over_knockouts_formate} depicts the individual cost factors for formate production. 
The biggest cost driver is substrate cost, amounting to 50 to 81\% of the total cost. 
For these small production capacities, 50\% is a typical value.
Further, the pH adjustment agent is a non-negligible factor, which matches the findings of \citep{Saur.2023} in their techno-economic analysis of itaconic acid production. 
Substrate costs are followed by investment costs, which consist almost entirely of reactor costs. 
The compressor investment cost, and also compressor power costs, are low due to the small reactor size.
Utilities are the lowest cost item, split roughly equally between power and cooling costs. 
This cost distribution is also observable for acetate and succinate production (see Figure~\ref{fig:cost_increasing_KOs_appendix}).
The only exception is succinate production with four knockouts, where reactor costs (\SI{1.0}{\USD\per\kg}) are larger than substrate costs (\SI{0.86}{\USD\per\kg}). 
In this case, the optimal solution implies three parallel reactors.
 
If one compares the results with the embedded core network with those with the embedded GEM, the shares are different.
For example, the pH adjusting agent is not a cost factor with the embedded GEM, in contrast to results with the embedded core network. 
Presumably, the GEM allows for intracellular pH adjustment. 
Another example of different cost allocations is formate production with one knockout (which is ATP synthase with GEM and also with core network embedded). 
The predicted process with the embedded core network operates with six reactors, \SI{0.002}{\cubic\meter} each, an oxygen mole flow of \SI{0.22}{\mol\per\hour} (that is, \SI{0.15}{\mol\per\liter\per\hour}), a specific power input of \SI{8.2}{\kilo\watt\per\cubic\meter} and a pH of 6.0. 
On the other hand, the process with the embedded GEM operates with one reactor of \SI{0.05}{\cubic\meter}, an oxygen mole flow of \SI{5.9}{\mol\per\hour} (that is, \SI{0.18}{\mol\per\liter\per\hour}), a specific power input of \SI{5.9}{\kilo\watt\per\cubic\meter} and a pH of 6.5. 
This example shows that the optimal operation point may differ depending on the chosen metabolic network. 
Moreover, the example depicts that the optimization may choose several small parallel reactors over one bigger reactor. 
The solver only returns one set of variables for the optimal objective value. 
However, other solutions may exist. 
For example, if the metabolic network contains two pathways that are equally effective, the solver arbitrarily picks one. 
Analogously, it is possible that a solution with fewer, larger reactors exists which is just not returned by the solver. 
Another reason could be that smaller reactors actually achieve lower costs by the specific power input, which is bigger with smaller culture volumes. 
The specific power input directly influences the oxygen uptake (see \eqref{eq:vantRiet}), which in turn is decisive for growth and product formation.

\subsection[Variation of production capacity]{Economies of scale come into play when the production capacity is increased}

\begin{figure}[h!tb]
     \centering
     \begin{subfigure}[b]{0.49\textwidth}
     \centering         \includegraphics[width=\textwidth]{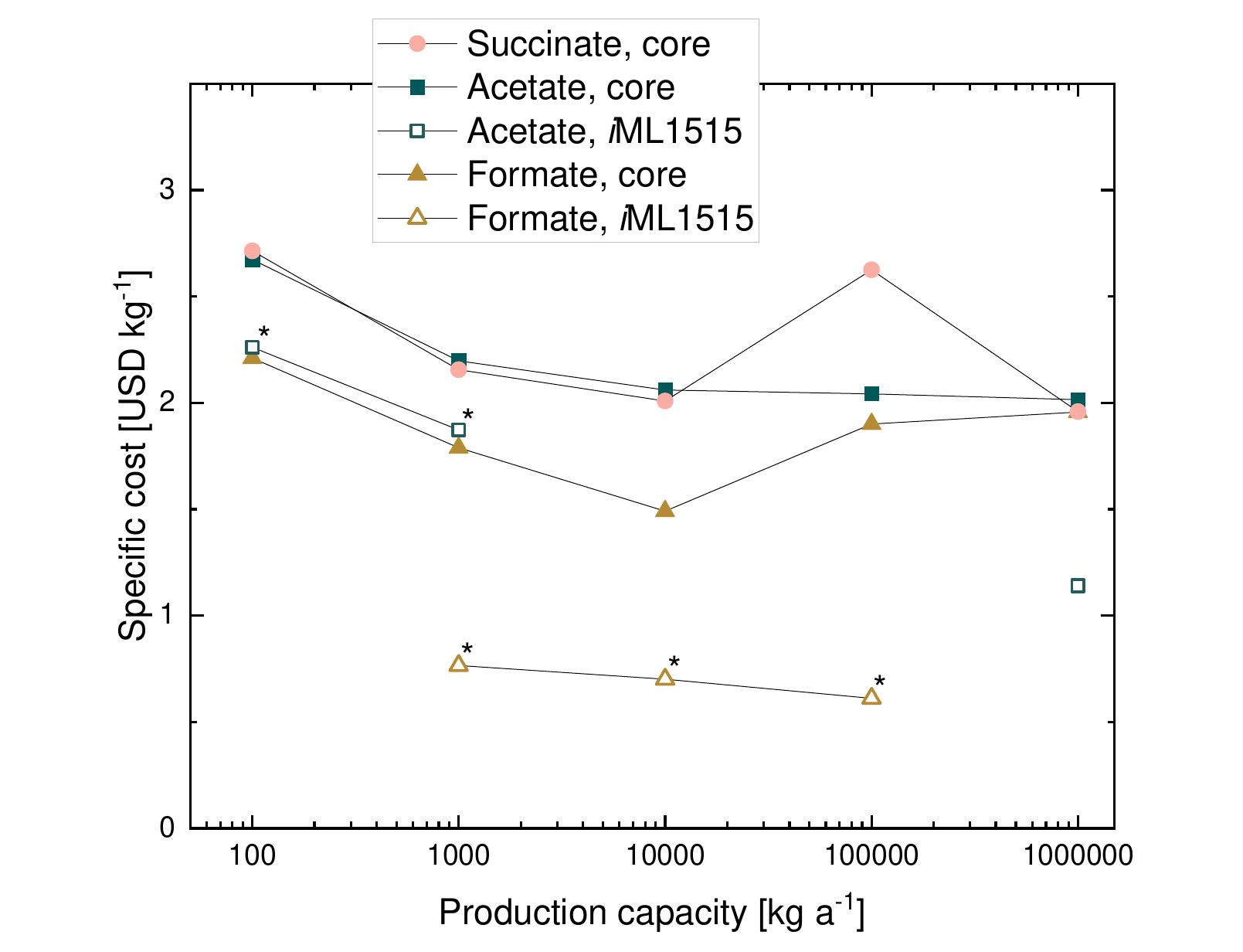}
         \caption{Total cost}
        \label{subfig:cost_over_capacity_total}
     \end{subfigure}
     \hfill
     \begin{subfigure}[b]{0.49\textwidth}
     \centering
     \includegraphics[width=\textwidth]{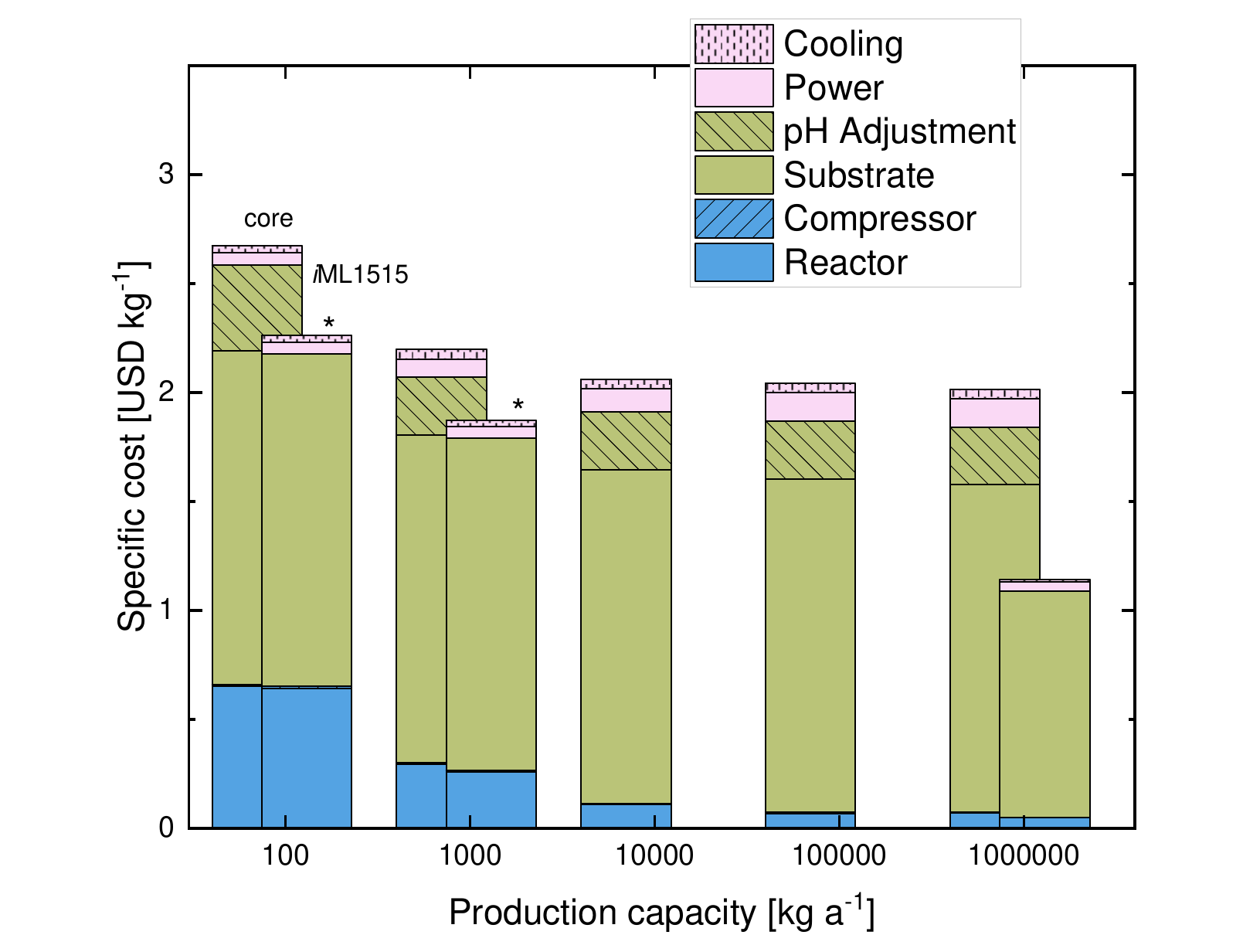}
     \caption{Acetate: cost factors}
       \label{subfig:cost_over_capacity_acetate}
      \end{subfigure}
     \caption{Influence of the production capacity on the specific production cost. The maximum allowable number of  knockouts is set to 3. The asterisk * indicates locally optimal solutions; all other solution points are globally optimal. The embedded networks are \textit{E. coli} core \citep{Orth.2010a} (indicated or bars in the back) and \textit{i}ML1515 \citep{Monk.2017} (indicated or bars in the front). USD: US-Dollar. (a) Total cost for acetate, succinate, and formate, (b) Acetate production cost divided into the individual cost factors. Blue, green, and pink color denotes investment cost, raw materials, and utilities, respectively.}
     \label{fig:cost_increasing_capacity}
 \end{figure}

In the next case study, we varied the production capacity from 100 to \SI{e6}{\kg\per\year} and examined the impact of this variation on the total specific cost and the individual cost factors. 
Figure~\ref{subfig:cost_over_capacity_total} depicts the total specific cost for succinate, acetate, and formate production.
With the embedded GEM, no global or local solution was found for succinate production.
Furthermore, the local solution for formate production of \SI{100}{\kg\per\year} was higher than the result with the embedded core network and therefore excluded from the figure. 

Economies of scale apply to acetate production, that is, the larger the production capacity, the smaller the specific cost. 
With the embedded core network, total cost decreases by 25\% from 100 to \SI{e6}{\kg\per\year}; with the embedded GEM, the decrease is 50\% from 100 to \SI{e6}{\kg\per\year}. 
Especially the reactor investment cost decreases with increasing production capacity (see Figure~\ref{subfig:cost_over_capacity_acetate}). 
Moreover, again, no pH adjustment is predicted with the embedded GEM. 
The knockout strategy with the embedded core network remains the same across the entire range of production capacity. 
As shown in Section~\ref{sec:cost_over_knockouts}, two knockouts are sufficient to reach a cost plateau for acetate production, namely ATP synthase and pyruvate kinase. 
Hence, the third chosen knockout has little effect. 
For \SI{e6}{\kg\per\year} acetate production with the embedded GEM, the suggested knockouts are also ATP synthase and pyruvate kinase. 
The other cases with the embedded GEM predict other knockout targets next to ATP synthase, however, these knockouts have little effect, as argued in Section~\ref{sec:cost_over_knockouts}.
Overall, we concluded that economies of scale apply. 
However, they only apply to investment costs (see \eqref{eq:invest_reactor} and \eqref{eq:invest_compressor}).
With substrate cost being the largest cost factor, however, they do not change the orders of magnitude. 
Considering economies of scale also in the procurement of raw materials might change the picture. 

The cases of formate and succinate production are more complex (see Figure~\ref{fig:cost_increasing_capacity_appendix}). 
With succinate, the production of \SI{e5}{\kg\per\year} marks an outlier with a change in knockout strategy. 
This case suggests deleting the reactions  6-phosphogluconolactonase, fumarase, and pyruvate dehydrogenase, whereas all other cases suggest removing acetate reversible transport and ATP synthase.
This change in knockout strategy leads to increased growth and decreased product flux, which, in turn, increases substrate and pH adjustment costs. 
Such behavior can occur when variable bounds become active in the optimization which were not active in the other cases.
The variable bounds may represent physically meaningful limits for the variables. 
We therefore examined the active bounds of this case.
The substrate flux is at its upper bound of \SI{10}{\milli\mol\per\gramdryweight\per\hour}; however, this is the fact in all cases. 
Further, the OUR is almost at its upper bound of \SI{250}{\milli\mol\per\liter\per\hour}.
However, this bound is also active in other cases, for example, \SI{e4}{\kg\per\year} acetate production, which behaves regularly. 
In summary, we could not explain the irregularity of this case.

With formate, the specific production cost decreases with increasing production capacity. 
With the embedded core network, however, they increase again at \SI{e5}{\kg\per\year} and \SI{e6}{\kg\per\year}.
This increase is associated with increased raw material costs; the investment costs continue to decline. 
With the embedded GEM, this increase in production cost is not evident. 
Hence, presumably, large production capacities are suboptimal for production with the \textit{E. coli} core network. 


\subsection[Simultaneous vs. sequential approach]{The simultaneous approach outperforms a sequential approach}

Our scientific hypothesis is that simultaneous optimization is advantageous over sequential optimization for microbe and bioreactor design. 
In order to examine the validity of this hypothesis, we defined a sequential approach with two steps. 
In the first step, OptKnock \citep{Burgard.2003} predicts knockouts for a specified target product and a fixed maximum allowable number of knockouts. 
The predicted knockouts are handed over to a second step. 
In the second step, SimulKnockReactor minimizes the reactor cost for the specified target product with fixed knockouts.
Note that already this approach is itself a novelty. 
This two-step procedure reflects the established approach to first optimize the microorganism and then design the bioreactor. 

We applied this sequential approach to formate, acetate, and succinate production with an embedded \textit{E. coli} core network \citep{Orth.2010a}. 
We chose to work with the core network to gain a full overview of the cases and have less impact by numerical issues. 
The production capacity was set to \SI{e3}{\kg\per\year} to make the results directly comparable with the results from Section~\ref{sec:cost_over_knockouts} (with embedded core network).

When performing the sequential approach, different outcomes may arise for the different target products and number of knockouts.
These outcomes are based on the feasibility of an optimization program. 
In general, an optimization program has a feasible solution if all equalities and inequalities are satisfied by one set of variables \citep{Edgar.2001}.
In coherence with this definition, the solver returns a statement of whether the optimization program has a feasible solution. 
Due to numerical reasons, however, this output may not reflect the mathematical truth; for example, the program may be infeasible even if the solver returns a solution. 
The case that the solver declares a program infeasible even if it is actually feasible, however, is unlikely in global optimization. 
Therefore, we reported feasibility based on the solver output. 

\begin{figure}[h!tb]
    \centering
    \includegraphics[width=0.75\linewidth]{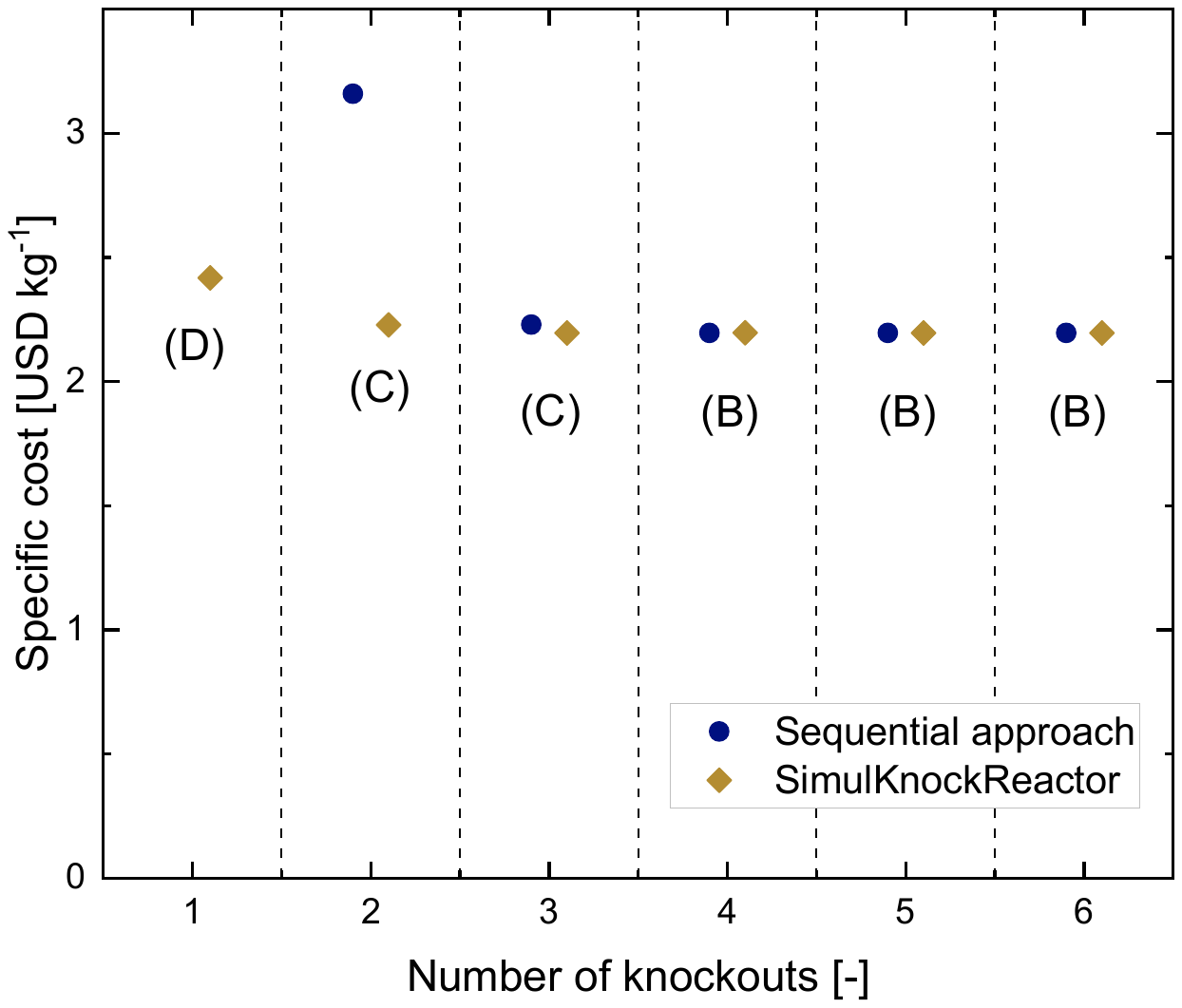}
    \caption{Specific cost of acetate over the number of knockouts for the sequential approach and SimulKnockReactor. Production capacity is set to \SI{e3}{\kg\per\year}; the embedded network is \textit{E. coli} core \citep{Orth.2010a}. Outcomes B, C, and D are defined in Table~\ref{tab:Cases_seq}. USD: US-Dollar.}
    \label{fig:sim-vs-seq_Acetate}
\end{figure}

In the different case studies, the solver may return that the sequential optimization program is infeasible in the first stage.
However, we did not observe this outcome in our case studies.
Next, the sequential optimization program can turn infeasible in the second stage. 
When this is the case, it is interesting to compare with the behavior using SimulKnockReactor. 
If the solver returns a feasible solution using SimulKnockReactor, this indicates that the proposed knockouts using OptKnock are not able to fulfill the production capacity requirements of the reactor stage, while the knockouts using SimulKnockReactor are.
The case where the optimization program using the sequential approach has a feasible solution and the optimization program using the simultaneous approach is infeasible is mathematically not possible and numerically very unlikely because the solver is free to choose the same knockouts as OptKnock also with the simultaneous approach. 
When the solver returns a feasible solution with both approaches, we compared the proposed knockout strategy. 
From the same proposed knockouts, it follows directly that the objective values are the same, too. 
With different proposed knockouts, the objective value with SimulKnockReactor can be equally good or better than with the sequential approach. 
The case that the objective values with SimulKnockReactor are worse is mathematically not possible since the same knockouts as OptKnock can be chosen. 
Nevertheless, it can still happen with the solver because the simultaneous optimization program is numerically very challenging.
The outcomes that we observed in our case studies are depicted in  Table~\ref{tab:Cases_seq}).
Exemplarily, Figure~\ref{fig:sim-vs-seq_Acetate} shows the specific cost of acetate production; the data for formate and succinate can be found in the Supplementary Information. 

\begin{table}[h!tb]
    \centering
    \caption{Observed outcomes based on the solver output for formate, acetate, and succinate production from one to six knockouts and a production capacity of \SI{e3}{\kg\per\year} with the embedded \textit{E. coli} core \citep{Orth.2010a} network.}
    \label{tab:Cases_seq}
    \resizebox{\textwidth}{!}{
    \begin{tabular}{c|c|c}
        \backslashbox{Sequential}{SimulKnockReactor} & feasible  & infeasible \\
        \hline
        \multirow{3}{*}{feasible}    & (A) same knockouts   & \multirow{3}{*}{not observed} \\
         & (B) different knockouts, same objective value & \\
         & (C) different knockouts, SimulKnockReactor better & \\
         \hline
        infeasible (2\textsuperscript{nd} stage) & (D) & (E) \\
    \end{tabular}
    }
\end{table}

The solver returns feasible solutions in more cases with SimulKnockReactor than with the sequential approach (outcome D).
Example cases for this outcome are acetate production with one knockout (see Figure~\ref{fig:sim-vs-seq_Acetate}), succinate production with more than three knockouts (below three knockouts, also SimulKnockReactor is infeasible), and formate production with one and two knockouts. 
Hence, the optimization program of the sequential approach appears to be more prone to infeasibility when few knockouts are allowed.

For the cases where the solver returns feasible solutions with both approaches, SimulKnockReactor objective values are equally good or better than the ones of the sequential approach (outcomes A to C).
The same knockouts (outcome A) are predicted for formate production with three knockouts. 
In tendency, the higher the maximum allowable number of knockouts, the more the objective values resemble each other, even if the predicted knockouts differ (outcome C). 
This can be explained by predicted knockouts that do not further improve the result, as is the case with acetate production with four to six knockouts (see Figure~\ref{fig:sim-vs-seq_Acetate}). 
Furthermore, with more knockouts allowed, the number of possible combinations that lead to the same objective value may increase. 

The biggest cost item is substrate cost, with SimulKnockReactor and with the sequential approach. 
The substrate flux, in turn, directly influences the product flux, which is the objective function of OptKnock.
Hence, currently, there is a strong correlation between the simultaneous and the sequential approach. 
This correlation is expected to be less strong if other cost drivers become more important, for example, through raised electricity or compressor costs.

\subsection[Biochemical vs. thermochemical formate production]{Industrial, thermochemical formate production is still cheaper but of the same order of scale as biochemical production}

\begin{figure}[h!tb]
    \centering
     \begin{subfigure}[b]{0.49\textwidth}
     \centering         
        \includegraphics[width=\textwidth]{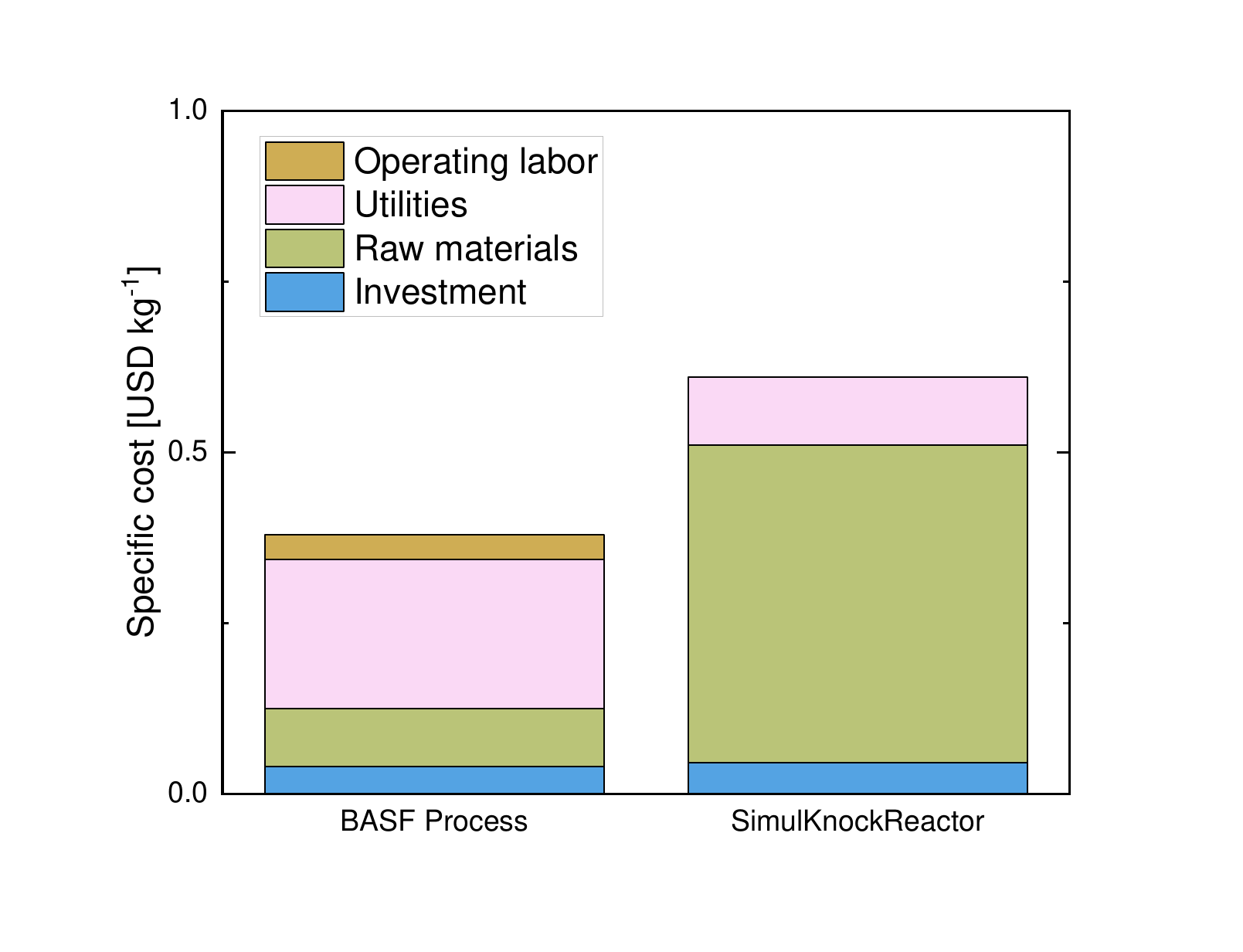}
         \caption{Cost factors}
    \label{subfig:cost_BASF}
     \end{subfigure}
     \hfill
     \begin{subfigure}[b]{0.49\textwidth}
     \centering
        \includegraphics[width=\textwidth]{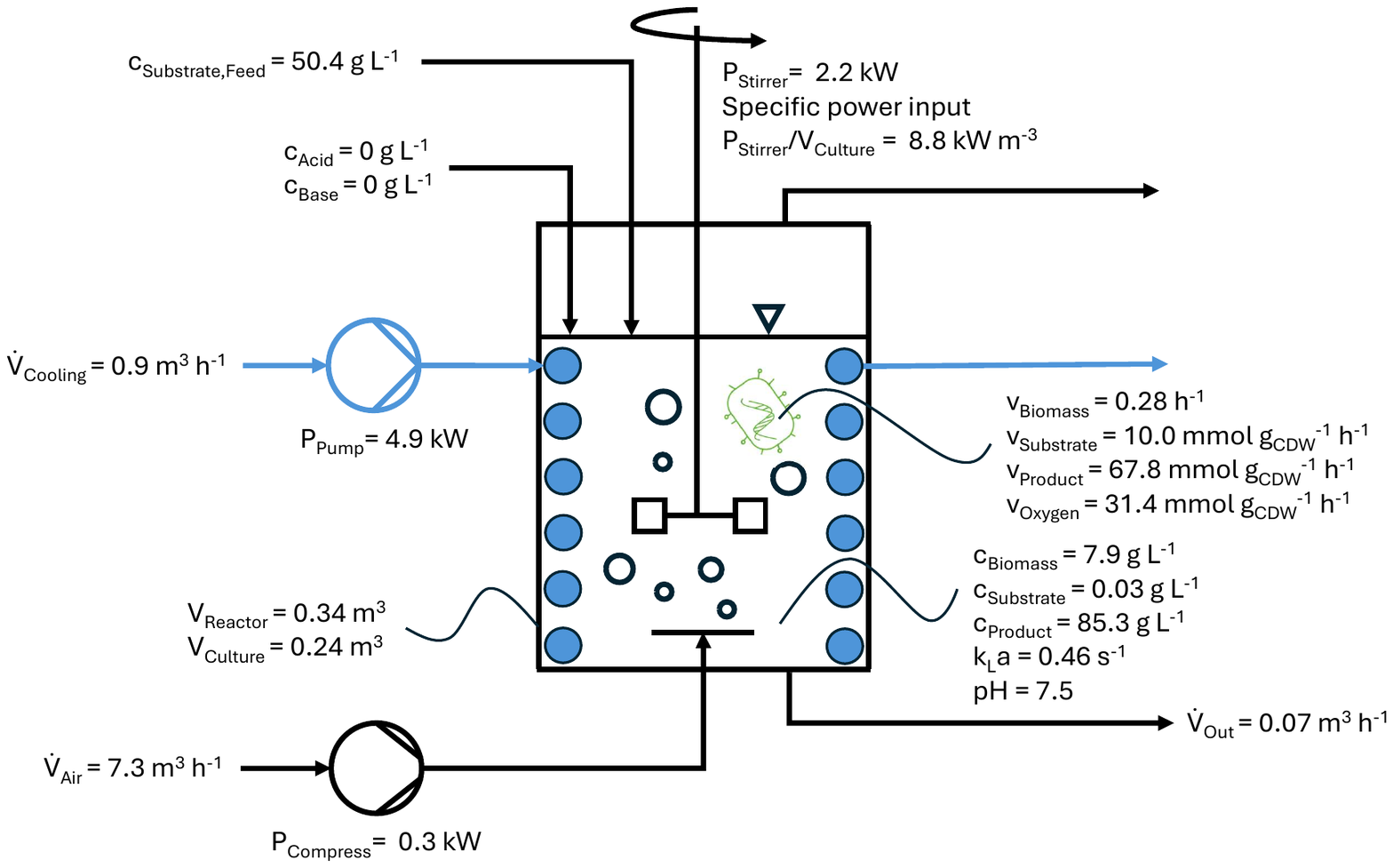}
     \caption{Predicted bioreactor and microbe design}
     \label{subfig:flowsheet_BASF}
      \end{subfigure}
    \caption{Comparison of an optimized thermochemical production process from BASF \citep{DaCunha.2018, Hietala.2003} and the biochemical production predicted with SimulKnockReactor. The production capacity is \SI{2.71e7}{\kg\per\year}. SimulKnockReactor is run for three maximum allowable knockouts with the embedded genome-scale metabolic model \textit{i}ML1515 \citep{Monk.2017}. USD: US-Dollar, c: concentration, P: power, V: volume. (a) Distribution of formate production costs (b) Bioreactor and microbe design variables predicted using SimulKnockReactor}
    \label{fig:BarChart_BASF_process}
\end{figure}

In our last case study, we evaluated the production cost of formate with a production capacity of \SI{2.71e7}{\kg\per\year}.
We, thereby, aimed to compare the cost of the biochemical production of formate with the cost of an established thermochemical production route. 
Moreover, we wanted to evaluate the prediction ability of SimulKnockReactor in a real-world example.
With a worldwide market share of 32.1\% in 2013 and a company-specific production process \citep{Hietala.2003}, the company BASF was chosen as a benchmark for thermochemical formate production.  
Production costs of the thermochemical formate production were taken from \citet{DaCunha.2018}, who optimized the total annual cost of the formic acid process by the BASF company \citep{Hietala.2003} in a bi-objective optimization. 
Production costs of the biochemical route were predicted with SimulKnockReactor using local optimization. 
The production capacity was set to \SI{2.71e7}{\kg\per\year}, equal to \citet{DaCunha.2018}, and the number of possible knockouts was set to three. 
The prediction was performed using the \textit{i}ML1515 network \citep{Monk.2017}.
In SimulKnockReactor, the reactors and the compressor amount to investment cost; substrate and pH adjusting agent sum up to raw materials; and power and cooling make up for utilities. 
Operating labor costs were not considered in SimulKnockReactor and are therefore not available. 
Figure~\ref{subfig:cost_BASF} depicts the production cost of formate in the BASF process and predicted by SimulKnockReactor. 

Biochemical production is 60\% more expensive than the thermochemical production.  
The biggest difference lies in the cost of raw materials. 
With SimulKnockReactor, the raw materials account for 76\% of the total cost, which is on the upper end of a typical range with large production capacities. 
However, the value highly depends on the process and the product.
In the large-scale production of itaconic acid, for example, \citet{Nieto.2020} found that the raw materials only account for around 12\% of total OPEX per year.
While the BASF process is based on carbon monoxide (\SI{0.07}{\USD\per\kg}), methanol (\SI{0.672}{\USD\per\kg}), and water (\SI{0.001}{\USD\per\kg}) \citep{DaCunha.2018}, the main substrate of the biochemical route is the more expensive glucose (\SI{0.79}{\USD\per\kg}) \citep{ChemAnalyst.2023_Glc}. 
Adjusting agents for pH are not needed and, thus, their costs are zero. 
The BASF process includes recycling streams, which increase the product yield. 
Recycling streams were not considered in the biochemical process (see Figure~\ref{subfig:flowsheet_BASF}).
The outlet stream only contains low substrate concentrations (\SI{0.03}{\gram\per\liter}), so substrate recycling is not of interest.  
However, recycling or retention of biomass could be interesting to decouple growth from the dilution rate and achieve lower residence times.
Moreover, an increased yield would be more favorable. 
Currently, the yield of formate per substrate on a carbon molar basis is 56\%.
However, the main side product is biomass, which is also needed for production.  

The utility costs are the largest cost factor in the BASF process, presumably due to the high pressure and temperature levels required for the thermochemical reaction and separation. 
This finding highlights the advantage of biochemical processes, which run at moderate temperature and pressure levels. 
However, it should be noted that the BASF process includes separation units and heat integration, which are not considered in SimulKnockReactor. 
This mismatch makes a deepened comparison difficult. 

The investment costs are similar in both processes. 
The biggest share with SimulKnockReactor is the reactor investment costs. 
Namely, in order to reach the production capacity, 545 reactors operate in parallel, with a volume of \SI{0.34}{\cubic\meter} each. 
As already discussed in Section~\ref{sec:cost_over_knockouts}, the more realistic alternative of fewer reactors with larger volumes each is presumably not chosen by the solver due to higher specific power input in small reactors.
Adaptation of the model by penalizing parallelization of reactors could potentially tackle this issue.  
This large amount of reactors is needed due to the large residence time, that is, \SI{3.5}{\hour}.
As a reminder, due to continuous fermentation conditions, the residence time is the inverse of the growth rate.
Decreasing the residence time would, hence, only be possible by increasing the growth rate. 
In contrast, in the BASF process, the residence time in reactor 1 and 2 are \SI{0.717}{\hour} and \SI{0.202}{\hour}, respectively.
Thus, the separation units presumably account for the largest share of investment cost. 

Lastly, despite the mismatch discussed above, we want to point out that the costs of the BASF process and the biochemical process predicted with SimulKnockReactor range in the same order of scale, which demonstrates the applicability of SimulKnockReactor as an early-stage process evaluation tool.  

\section{Conclusion and outlook}

We presented a bilevel optimization formulation, named SimulKnockReactor, which combines bioreactor design with microbial strain design. 
It is an extension of our previous work SimulKnock \citep{Ziegler.2024}.
The upper level minimizes total bioreactor costs and includes the mass and energy balances of the fermentation, as well as the dimensioning of the reactor, the aeration, the cooling, and the pH control. 
We referred to this level as bioreactor design, but excluded the selection of specific equipment, such as the type of stirrer.
The lower level maximizes the growth of the microorganism and includes metabolite mass balances based on the metabolic model of the organism. 
We considered gene deletion as a genetic modification of the microbial strain. 
The bilevel optimization formulation was reformulated to a single-level program resulting in a mixed-integer nonlinear optimization program. 
It was solved using global and local optimization techniques. 

Our results of the exemplary case studies with the \textit{E. coli} core network \citep{Orth.2010a} and the GEM \textit{i}ML1515 \citep{Monk.2017} showed that substrate cost was the biggest cost factor.
Further, we observed that the extent of pH adjustment cost depended on the embedded network. 
With an increasing number of knockouts, the production cost decreased until it reached a plateau. 
As production capacity increased, economies of scale occurred, especially in the investment costs for reactors.  
In comparison with a sequential approach, where OptKnock \citep{Burgard.2003} was applied and then the reactor cost was minimized, the solver returned feasible solutions in more case studies with SimulKnockReactor.
This indicated that, with the knockouts predicted using OptKnock, the microorganism could not fulfill the requirements of the reactor stage. 
Further, the total reactor cost with SimulKnockReactor was equal to or lower than with the sequential approach, which highlights the techno-economic advantages of this formulation and its predicted knockout strategies.  

We proved the solvability of SimulKnockReactor also for industrially relevant production capacities. 
However, the amount of reactors predicted by the solver was unrealistically high, presumably due to better oxygen distribution in small reactors. 
Additional constraints penalizing parallel operation could tackle this issue. 
To improve the accuracy of our aeration model, the employed van't Riet correlation to calculate the $k_La$-value should be supported with experimental values or a strain-specific correlation could be employed. 
Furthermore, the hydrostatic pressure in the tank could be introduced as a variable. 
Taking into account the effects of aeration on agitation power input would further refine the modeling. 
The accuracy of the pH model could be increased by incorporating metabolic models including pH \citep{Du.2019}, or by rigorous modeling of the pH, as suggested by \citet{Walz.2017} and \citet{Bree.2020}. 
Moreover, our assumption that genetic modifications do not influence the impact of pH on growth should be experimentally validated.  
The modeling of the cooling could be extended by considering the heat transfer surface of the cooling coils and a sensible upper bound on the mass flow of the cooling agent or pumping power.

Depending on the size of the metabolic network, the runtimes ranged from a few minutes using up to 8 threads (for core network) to several hours on 48~threads, with difficulties of convergence within the set time limit of \SI{10}{\hour}(for GEM). 
To decrease computation times and tackle convergence issues, machine learning could be included to substitute the GEM, to include experimental data (e.g., on pH influence), or to learn nonlinear functions. 

Beyond that, many, possibly modular extensions are conceivable.
An extension could be the selection of specific devices for the bioreactor and the surrounding equipment in order to do full justice to the concept of bioreactor design. 
The inclusion of product inhibition, e.g., using a linear correlation of maximal growth with product concentration  \citep{Straathof.2023}, could be interesting.
Furthermore, more aspects could be included, e.g., considerations about upstream and downstream units \citep{Konstantinos.2022, TafurRangel.2022, Ploch.2019, CortesPena.2020}, as well as marketing and administration, transportation, and research and development, as indicated by \citep{Doran.1995, Liu.2020}.
Using multi-objective optimization, the global warming impact could be considered parallel to cost estimation, similar to the work of \citet{Konig.2020}.
The inclusion of co-feeding, up-regulation, and down-regulation of genes \citep{Pharkya.2006}, as well as dynamic operation \citep{Mahadevan.2002}, i.e., batch and fed-batch, stay interesting.
With increasing available computational power and similar to \citet{Klamt.2004}, the lower level could be replaced by elementary mode analysis \citep{Schuster.1994} that describes the microbe without specifying an artificial cellular objective. 

Lastly, we want to emphasize again that, given the few publications that exist on model-based bioreactor design in general and bioreactor design coupled with the cellular level in particular, the presented formulation SimulKnockReactor can also enhance these fields of research.  
For example, it should be possible to run SimulKnockReactor without considering genetic modifications if the production of the target chemical is already growth-coupled.
In this use case, the number of allowed knockouts would be set to zero.
Presumably, the optimization of the reactor conditions will activate one path versus another at the cellular level, and, thus, refine the bioreactor optimization.
For a quick, first approximation of bioreactor cost and size, the substrate concentration in the reactor could be set to zero such that pH and kinetics are not considered. 
Even if all cellular variables at the upper level were replaced again by common correlations--- making the lower level obsolete---SimulKnockReactor could contribute to establishing the term and field of optimal bioreactor design.

\section*{Statements}
\subsection*{Author contributions}

CRediT: \\
ALZ: Conceptualization, Methodology, Software, Validation, 
Investigation, Writing - Original Draft, Visualization\\
MDS: Methodology, Software, Validation, Investigation, Writing - Review and Editing\\
TP: Methodology, Software, Investigation, Visualization, Writing - Review and Editing\\
TS: Methodology, Writing - Review and Editing\\
JM: Methodology, Writing - Review and Editing, Funding acquisition \\
AM: Conceptualization, Methodology, Resources, Writing- Review and Editing, Supervision, Funding acquisition

\subsection*{Acknowledgments} 
Computations were performed with computing resources granted by RWTH Aachen University under project thes1557.
We thank Luis Schrade and Ilias Ouzlim for their preliminary work on setting up the optimization formulation.
We thank Clemens Kortmann for making the code compatible with the code framework of SimulKnock.

\subsection*{Funding}
This project was funded by the Deutsche Forschungsgemeinschaft (DFG, German Research Foundation) under Germany´s Excellence Strategy – Cluster of Excellence 2186 ``The Fuel Science Center'' – ID: 390919832.

\subsection*{Declaration of conflicts}
    No conflicts to declare.

\subsection*{Data availability statement}
The implementation of the optimization formulation, the data pre-processing, and the interface with the solver Gurobi are openly available in our GitLab repository at \url{https://git.rwth-aachen.de/avt-svt/public/simulknock}.

\section*{Nomenclature}

\subsection*{Abbreviations}
\begin{longtable}[l]{p{45pt} p{250pt} }
	conc. & concentration \\
    CPM & cardinal pH model \\
    CSTR & continuous stirred-tank reactor \\
    FBA	 & flux balance analysis \\
    GEM  & genome-scale metabolic model\\
    OTR & oxygen transfer rate \\
    OUR & oxygen uptake rate \\
    substr. & substrate \\
    USD & US-Dollar\\
\end{longtable}

\subsection*{Symbols}
\begin{longtable}[l]{p{45pt} p{250pt} }
    $A$ & surface [\si{\square\meter}] \\
    $\boldsymbol{B}$ & mapping matrix [-]\\
    $c$ & concentration [\si{\gram\per\liter}]\\
    $c_{bio}$ & biomass concentration [\si{\gramdryweight\per\liter}]\\
    $c^{mol}$ & molar concentration [\si{\mol\per\liter}]\\
    $c_{p,heat}$ & specific heat capacity [\si{\joule\per\kg\kelvin}] \\
    $C$ & specific cost [USD \si{\per\kilogram} or USD \si{\per\cubic\meter} or USD  kWh$^{-1}$ or USD \si{\per\hour} or USD] \\
    $D$ & diameter [\si{\meter}] \\
    $f$ & fixed fraction [-]\\
    $g$ & gravity (9.81) [\si{\meter\per\square\second}] \\
    $H$ & height [\si{\meter}] \\
    $\Delta h_{comb}$ & enthalpy of combustion [\si{\kilo\joule\per\gram}] \\
    $k$ & heat transfer coefficient [\si{\watt\per\square\meter\per\kelvin}] \\
    $k_La$ & volumetric oxygen mass transfer coefficient [\si{\per\second}] \\
    $K$ & dissociation constant [\si{\mol\per\liter}] \\
    $K_{S}$ & Monod constant [\si{\gram\per\liter}] \\
    $m$ & mass [\si{kg}] \\
    $\dot{m}$ & mass flow [\si{kg\per\hour}] \\
    $M$ & molar mass [\si{\gram\per\milli\mol}]\\
    $n_R$ & number of reactors set up in parallel [-] \\
    $\dot{n}$ & mole flow [\si{\mol\per\hour}]\\
    $OTR$ & oxygen transfer rate [\si{\gram\per\liter\per\hour}]\\
    $OUR$ & oxygen uptake rate [\si{\gram\per\liter\per\hour}] \\
    $p$ & pressure [\si{\pascal}]\\
    $P$ & power input [kW]\\
    $pH$ & pH [-] \\
    $q$ & relation height to diameter [-] \\
    $\dot{Q}$ & heat flow [\si{\kilo\watt}] \\
    $R$ & ideal gas constant (8.3145)[\si{\kilo\joule\per\kilo\mol\per\kelvin}]\\
    $\boldsymbol{S}$ & stoichiometric matrix [-] \\
    $t$ & operating time [\si{\hour}] \\
    $T$ & temperature [\si{\kelvin}]\\
    $u_s$ & superficial gas velocity [\si{\meter\per\second}] \\
	$v$	& flux [\si{\milli\mol\per\gramdryweight\per\hour}] \\
    $v_{bio}$ & growth rate [\si{\per\hour}] \\
    $V$ & volume [\si{\cubic\meter}]\\
    $\dot{V}$ & volume flow [\si{\cubic\meter\per\hour}]\\
    $y$ & knockout variable [-]
\end{longtable}

\subsection*{Greek letters}
\begin{longtable}[l]{p{45pt} p{400pt} }
    $\alpha$ & polytropic exponent [-] \\
    $\beta$ & parameter that correlates reactor volume with culture volume [-] \\
    $\gamma$ & growth inhibition factor [-]\\
    $\delta$ &parameter that correlates OUR with metabolic heat [\si{\kilo\watt\hour\per\mol}] \\
    $\epsilon$ & holdup [-] \\
    $\eta$ & efficiency number [-]\\
    $\kappa$ & maximum allowable number of knockouts [-]\\
    $\nu$ & stoichiometric coefficient [-] \\
    $\rho$ & density [\si{\kilogram\per\cubic\meter}]\\
    $\varphi$ & parameter that correlates oxygen concentration with oxygen solubility [-] \\
\end{longtable}

\subsection*{Subscripts}
\begin{longtable}[l]{p{45pt} p{250pt} }
    $bio$ & biomass\\
    $CDW$ & cell dry weight\\
    $comb$ & combustion \\
    $compress$ & compressor \\
    $cross$ & cross-sectional \\
    $equil$ & equilibrium \\
    $H^+$ & hydrogen ions \\
    $i$ & summation index \\
    $met$ & metabolic \\
    $O_2$ & oxygen \\
	$P$ & product\\
    $R$ & reactor\\
    $ref$ & reference \\
    $s$ & superficial \\
    $S$	 & substrate \\
    $WT$ & wild-type\\
\end{longtable}

\subsection*{Superscripts}
\begin{longtable}[l]{p{45pt} p{250pt} }
	$m$ & number of metabolites\\
    $max$ & maximum \\
    $min$ & minimum \\
    $mol$ & molar \\
	$n$	 & number of irreversible reactions \\
    $opt$ & optimal \\
	$r$ & number of (possibly reversible) reactions\\
\end{longtable}

\clearpage
\appendix
\setcounter{section}{0}
\renewcommand{\thesection}{A.\arabic{section}}
\setcounter{table}{0}
\renewcommand{\thetable}{A.\arabic{table}}
\setcounter{figure}{0}
\renewcommand{\thefigure}{A.\arabic{figure}}
\setcounter{equation}{0}
\renewcommand{\theequation}{A.\arabic{equation}}

\section{Tabulated metabolic and fermentation reaction parameters}
\begin{table}[h!tb]
    \centering
    \caption{Tabulated fermentation parameters, already presented and described in SimulKnock \citep{Ziegler.2024}. Ref.: reference}
    \resizebox{\textwidth}{!}{
    \begin{tabular}{lllll}
    \hline
         Symbol& Description& Value & Unit & Ref.   \\
    \hline
        $f$ & share of wild-type growth rate & 0.1 & - & - \\
        $K_S$ & Monod constant& 0.044&\si{\gram\per\liter}& \citep{Wick.2001} \\
        $v_{ATP}^{lower}$ & ATP maintenance flux threshold & 6.86 & \si{\milli\mol\per\gramdryweight\per\hour} & \citep{Monk.2017} \\
        $v_{bio}^{max}$ & maximum growth rate & 0.73 & \si{\per\hour} & \citep{Wick.2001}\\
        $v_S^{upper}$ & upper bound substrate uptake & 10 & \si{\milli\mol\per\gramdryweight\per\hour}  & -\\  
    \hline
    \end{tabular}
    }
    \label{tab:fermentation_parameter_values}
\end{table}

\section{Tabulated bioreactor parameters}

\begin{table}[h!tb]
    \centering
    \caption{Tabulated bioreactor parameters}
    \resizebox{\textwidth}{!}{
    \begin{tabular}{lllll}
    \hline
         Symbol& Description& Value & Unit & Reference   \\
    \hline
        $\alpha$ & polytropic exponent & 1.4 & - & \citep{Luyben.2018} \\
        $\beta$ & headspace parameter & 1.2 & - & \citep{VantRiet.1991} \\
        $\delta$ & metabolic heat correlation parameter & 0.13 & [\si{\kilo\watt\hour\per\mol}] & \citep{Doran.1995, Knoll.2005, Cooney.1969} \\
        $\Delta h_{comb,S}$ & enthalpy of combustion of glucose & 2805.7 & \si{\joule\per\milli\mol} & \citep{Doran.1995} \\
        $\Delta h_{comb,NH_3}$ & heat of combustion of ammonia & -382.6 & \si{\joule\per\milli\mol} & \citep{Doran.1995} \\
        $\Delta h_{biomass}$ & heat of combustion of biomass & -552 & \si{\joule\per\milli\mol} & \citep{Doran.1995} \\
        $\eta_{compress}$ & compressor efficiency & 0.7 & - & \citep{Cardoso.2020} \\
        $\eta_{pump}$ & pump efficiency & 0.7 & - & - \\
        $\rho_{cooling}$ &cooling water density & 1000 & \si{\kilogram\per\cubic\meter} &-  \\
        $\rho_{culture}$ & culture density & 1000 & \si{\kilogram\per\cubic\meter} & \citep{Mayer.2023, MartinezSalas.1981} \\
        $\rho_{O_2}$ & oxygen density at $T_{culture}$ & 1.27 & \si{\kilogram\per\cubic\meter} & \citep{Stephan.2019} \\
        $\varphi$ & correlation parameter dissolved oxygen & 0.3 & - & \citep{Bernard.2001} \\
        $c_{O_2,equil}$ & oxygyen solubility at $T_{culture}$ & 0.0068 & \si{\gram\per\liter} & \citep{Baburin.1981,Vendruscolo.2012} \\
        $c_{p,heat}$ & specific heat capacity cooling water at $T_{mean}$ & 4184 & \si{\joule\per\kg\per\kelvin}& \citep{Lemmon.2024} \\
        $C_{acid}$ & sulfuric acid specific cost & 0.064 & USD \si{\per\kg} & \citep{ChemAnalyst.2024_H2SO4}\\
        $C_{base}$ &sodium hydroxide specific cost & 0.29 & \si{\USD\per\kg} & \citep{ChemAnalyst.2023_NaOH} \\
        $C_{cooling}$ &specific cost cooling water & 0.17 & \si{\USD\per\cubic\meter} & \citep{Wang.2021} \\
        $C_{Power}$ & power specific cost & 0.06 & USD kWh$^{-1}$ & \citep{Humbird.2017} \\
        $C_{R,ref}$ & reference reactor cost & 40.000 & USD & \citep{Powell.2008} \\
        $C_S$ & glucose specific cost & 0.79 & \si{\USD\per\kg} & \citep{ChemAnalyst.2023_Glc} \\
        $g$ & gravity & 9.81 & \si{\meter\per\square\second} & - \\
        $M_{acid}$ & molar mass sulfuric acid & 98.1 & \si{\gram\per\mol} & - \\
        $M_{base}$ & molar mass sodium hydroxide & 40.0 & \si{\gram\per\mol} & - \\
        $M_{biomass}$ & molar mass biomass & 24.6 & \si{\gram\per\mol} & \citep{Doran.1995}\\
        $n_R^{upper}$& maximum number of reactors in parallel & 10 & - & - \\
        $p^0$ & atmospheric pressure & 1 & \si{\bar}& - \\
        $p_{air}$ & pressure of pressurized air & 2.5 & \si{\bar} & \citep{Luyben.2018}\\
        $pH^{lower}$ & lower bound pH & 4.5 & - & \citep{Du.2019} \\
        $pH^{max}$ & maximum pH & 9 & - & \citep{Rosso.1995} \\
        $pH^{min}$ & minimum pH & 4 & - & \citep{Rosso.1995} \\
        $pH^{opt}$ & optimum pH & 7 & - & \citep{Rosso.1995} \\
        $pH^{upper}$& upper bound pH & 7.5 & -& \citep{Akkermans.2017}\\
        $q$ & relation height to diameter & 2.5 & - &\citep{Humbird.2017} \\
        $R$ & ideal gas constant & 8.3145& \si{\kilo\joule\per\kilo\mole\per\kelvin}& -\\
        $t_{amortization}$& amortization period & 15 & a & \citep{Powell.2008} \\
        $t_{annual}$ & annual production hours & 8400 & h a$^{-1}$ & \citep{Humbird.2017}\\
        $T_{culture}$ & temperature fermentation culture & 303.15 &K & - \\
        $T_{cooling}^{in}$ & input temperature cooling water & 288.15 & K & - \\
        $T_{cooling}^{out}$ & output temperature cooling water & 298.15 & K & - \\
        $T_{mean}$ & mean temperature cooling water& 293.15 & K & - \\
        $v_{bio}^{opt}$& maximum growth rate & 0.73 & \si{\per\hour} & \citep{Wick.2001}\\
        $V_{R,ref}$ & reference reactor volume & 0.5 & \si{\cubic\meter} & \citep{Powell.2008}\\
        $V_R^{upper}$ & upper bound on the reactor size & 10 &\si{\cubic\meter} &- \\
    \hline
    \end{tabular}
    }
    \label{tab:bioreactor_parameter_values}
\end{table}

\clearpage
\section[Cost over knockouts]{Cost factors for acetate and succinate production over number of knockouts}

 \begin{figure}[h!tb]
     \centering
     \begin{subfigure}[b]{0.49\textwidth}
        \centering
        \includegraphics[width=\textwidth]{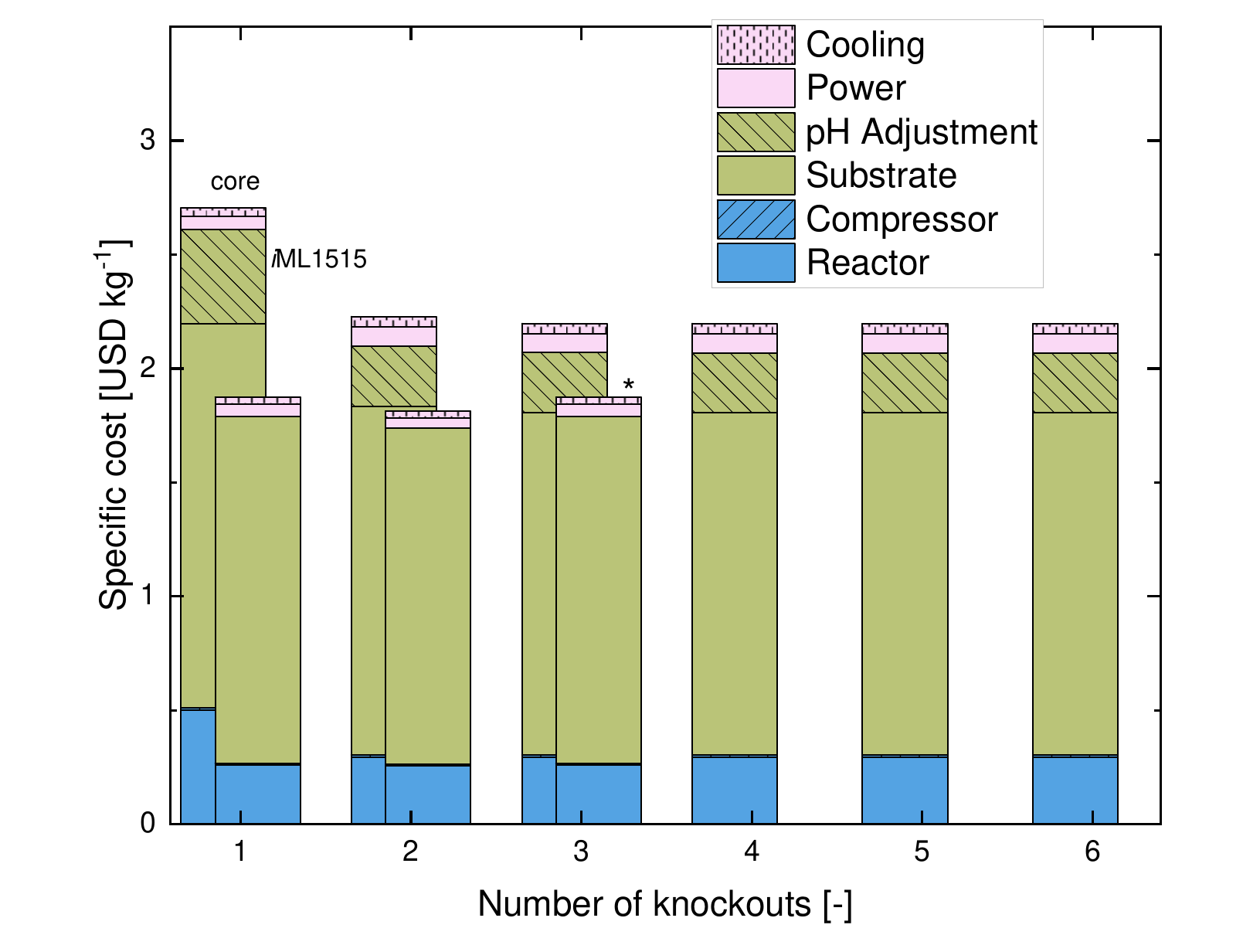}
        \caption{Acetate: cost factors}
        \label{subfig:cost_over_knockouts_acetate}
    \end{subfigure}
    \hfill
     \begin{subfigure}[b]{0.49\textwidth}
     \centering
         \includegraphics[width=\textwidth]{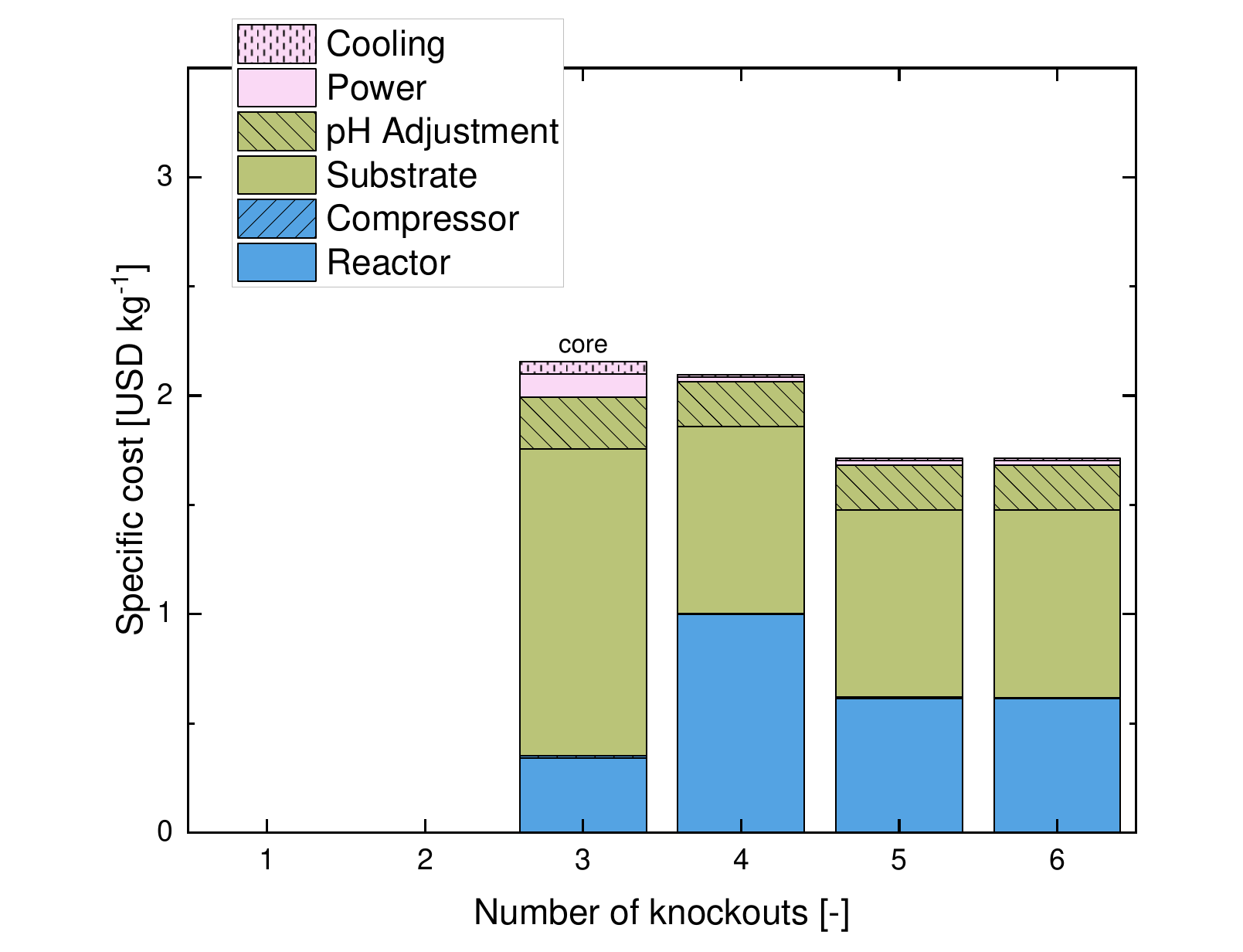}
         \caption{Succinate: cost factors}
         \label{subfig:cost_over_knockouts_succinate}
     \end{subfigure}
     \caption{Influence of the maximum allowable number of knockouts on the specific production cost. The production capacity was set to \SI{e3}{\kg\per\year}. The asterisk * indicates locally optimal solutions; all other solution points are globally optimal. The embedded networks are \textit{E. coli} core \citep{Orth.2010a} (indicated or bars in the back) and \textit{i}ML1515 \citep{Monk.2017} (indicated or bars in the front) USD: US-Dollar. (a) Acetate and (b) succinate production cost divided into the individual cost factors. Blue, green, and pink color denotes investment cost, raw materials, and utilities, respectively. }
     \label{fig:cost_increasing_KOs_appendix}
 \end{figure}

\section[Cost over capacity]{Cost factors for formate and succinate production over production capacity}

 \begin{figure}[h!tb]
     \centering
     \begin{subfigure}[b]{0.49\textwidth}
        \centering
        \includegraphics[width=\textwidth]{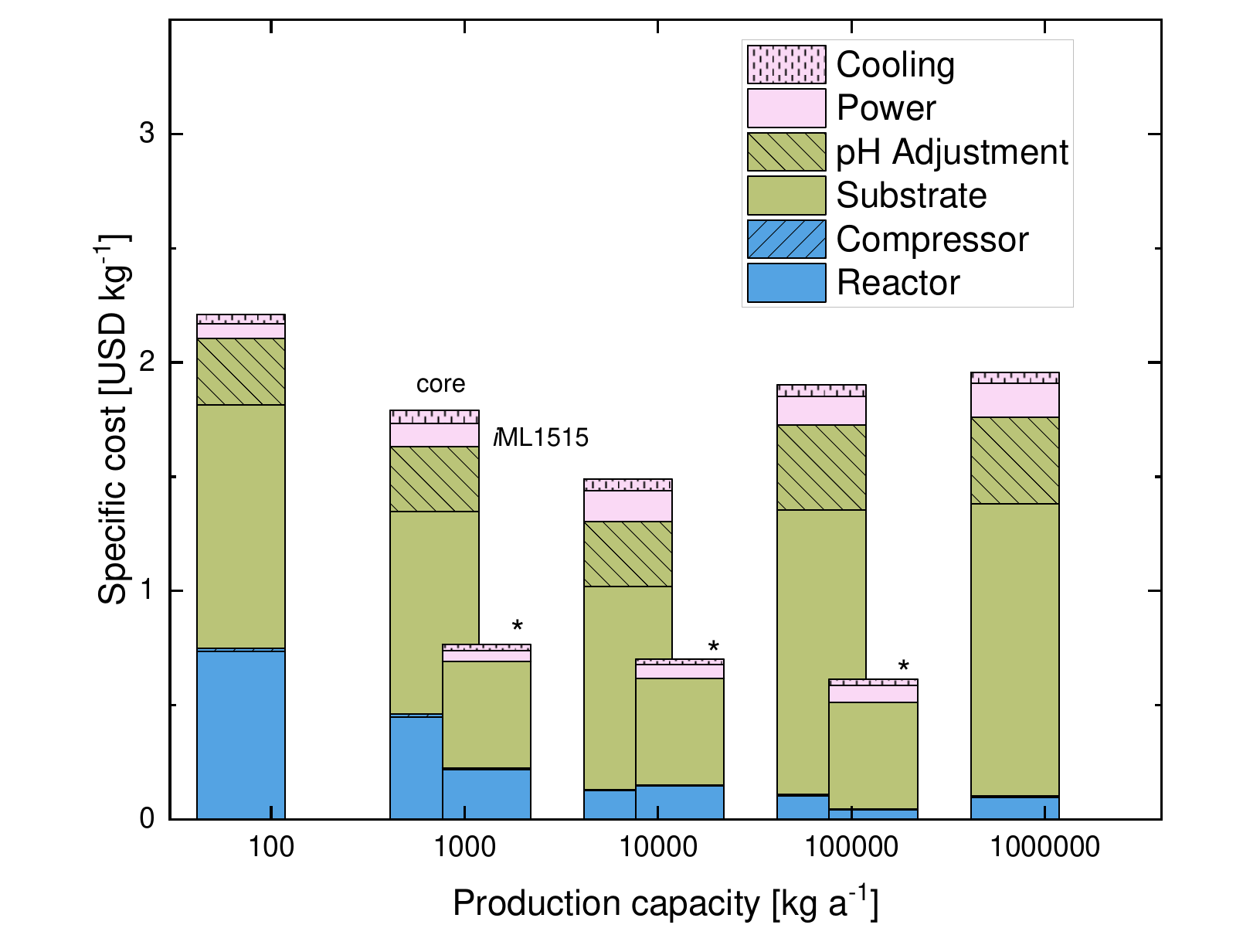}
        \caption{Formate: cost factors}
        \label{subfig:cost_over_capacity_formate}
    \end{subfigure}
    \hfill
     \begin{subfigure}[b]{0.49\textwidth}
     \centering
         \includegraphics[width=\textwidth]{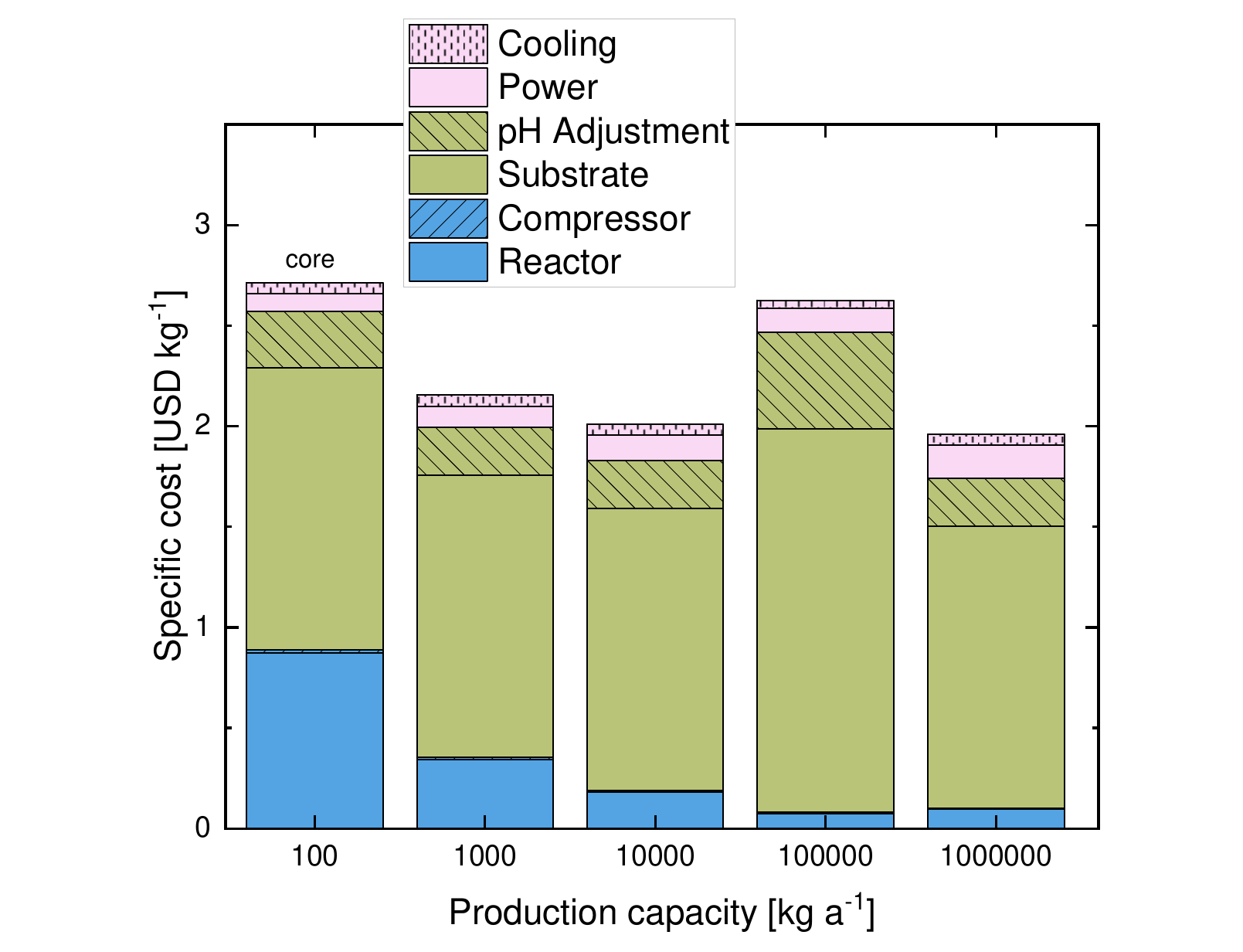}
         \caption{Succinate: cost factors}
         \label{subfig:cost_over_capacity_succinate}
     \end{subfigure}
     \caption{Influence of the production capacity on the specific production cost. The maximum allowable number of  knockouts was set to 3. The asterisk * indicates locally optimal solutions; all other solution points are globally optimal. The embedded networks are \textit{E. coli} core \citep{Orth.2010a} (indicated or bars in the back) and \textit{i}ML1515 \citep{Monk.2017} (indicated or bars in the front) USD: US-Dollar. (a) Formate and (b) succinate production cost divided into the individual cost factors. Blue, green, and pink color denotes investment cost, raw materials, and utilities, respectively. }
     \label{fig:cost_increasing_capacity_appendix}
 \end{figure}
\clearpage
\bibliographystyle{elsarticle-harv}
\bibliography{bibliography.bib}
\end{document}